\documentclass{osa-article}

\journal{oe}

\articletype{Research Article}

\usepackage{savesym}
\usepackage{wrapfig}
\usepackage{multirow}

\savesymbol{openbox}
\usepackage{amsthm}

\theoremstyle{definition}
\newtheorem{remark}{Remark}

\newtheorem{proposition}{Proposition}

\begin{document}

\title{A Local Fourier Slice Equation}

\author{Christian Lessig\authormark{1}}

\address{\authormark{1}Institut f{\"u}r Simulation und Graphik, Otto-von-Guericke-Universit{\"a}t Magdeburg, Germany}

\email{\authormark{*}lessig@isg.cs.uni-magdeburg.de} 


\begin{abstract}
We present a local Fourier slice equation that enables local and sparse projection of a signal.
Our result exploits that a slice in frequency space is an iso-parameter set in spherical coordinates.
Therefore, the projection of suitable wavelets defined separably in these coordinates can be computed analytically, yielding a sequence of wavelets closed under projection.
Our local Fourier slice equation then realizes projection as reconstruction with ``sliced'' wavelets with computational costs that scale linearly in the complexity of the projected signal.
We numerically evaluate the performance of our local Fourier slice equation for synthetic test data and tomographic reconstruction, demonstrating that locality and sparsity can significantly reduce computation times and memory requirements.
\end{abstract}


\section{Introduction}
\label{sec:introduction}

The Fourier slice theorem~\cite{Radon1917,Bracewell1956,Bracewell1990,Garces2011} plays an important role in many optical applications, for example medical imaging~\cite{Kak2001,Herman2009,Epstein2007}, plenoptic cameras~\cite{Ng2005},  radio astronomy~\cite{Bracewell1956,Bracewell1967}, and (electron) microscopy~\cite{Crowther1970,Levoy2006a}.
We introduce an analogue of the theorem that is localized in space and frequency and that, among other things, enables the local projection of a signal $f(\pmb{x})$ from a compressed representation.

\begin{wrapfigure}[15]{r}[0pt]{0.4\columnwidth}
  \centering
  \vspace{-11pt}
  \includegraphics[width=0.4\columnwidth]{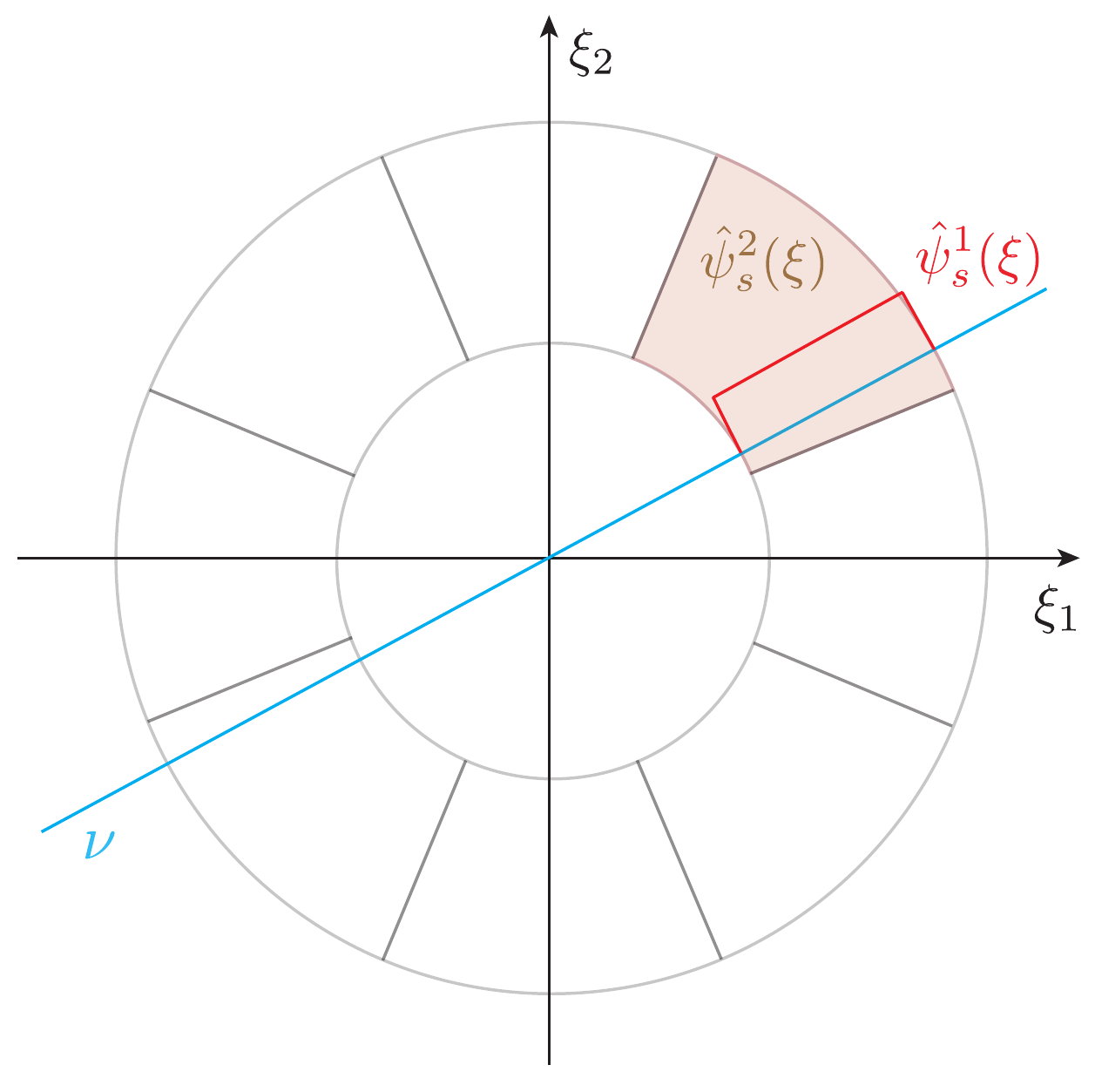}
  \vspace{-10pt}
  \caption{Conceptual view of our ansatz.}
  \label{fig:conceptual}
\end{wrapfigure}

Instead of the Fourier transform $\mathcal{F}(f) = \hat{f}(\pmb{\xi})$  used in the classical slice theorem,
\begin{align}
  \label{eq:fourier_slice}
  f_{\nu}(\pmb{y}) = \int_{\mathbb{R}_{\nu}} f(\pmb{x}) \, d\nu = \mathcal{F}_{P_{\nu}}^{-1} \big( \hat{f} \vert_{P_{\nu}} \big) ,
\end{align}
with $\nu$ being the direction along which the projection is performed and $P_{\nu}$ the (hyper-)plane orthogonal to it, our work relies on the polar wavelet representation of a signal,
\begin{align}
  \label{eq:wavelet_rep}
  f(\pmb{x}) = \sum_{s \in \mathcal{I}} f_s \, \psi_s^n( \pmb{x} ) ,
\end{align}
where $\psi_s^n( \pmb{x} )$ is a polar wavelet function in $\mathbb{R}_x^n$ whose Fourier transform is separable in polar coordinates.
Here $s = (j_s, k_s, t_s)$ is a multi-index that, in general, describes scale $j_s$, translation $k_s$, and orientation $t_s$.
Using that the restriction $\hat{f} \vert_{P_{\nu}}$ in Fourier space is along an iso-parameter set in polar coordinates, we show that polar wavelets form a sequence closed under projections,
\begin{align}
  \label{eq:psi_chain}
  \cdots \xrightarrow{\ \textrm{proj.} \ } \psi_{j,k,t}^n( \pmb{x} ) \xrightarrow{\ \textrm{proj.} \ } \psi_{j,k^{\nu},t^{\nu}}^{n-1}( \pmb{x} ) \xrightarrow{\ \textrm{proj.} \ } \cdots
\end{align}
``Slicing'' an $n$-dimensional polar wavelet $\psi_{j,k,t}^n(\pmb{x})$ thus yields an $(n-1)$-dimensional one $\smash{\psi_{j,k^{\nu},t^{\nu}}^{n-1}(\pmb{x})}$ at the same scale $j$ and with projected location $k^{\nu} = k - (k \! \cdot \! \nu) \nu$ and orientation $t^{\nu} = t - (t \! \cdot \! \nu) \nu$.
Furthermore, the wavelet $\smash{\psi_{j,k^{\nu},t^{\nu}}^{n-1}(\pmb{x})}$ does not depend on the projection direction, up to a scalar factor, and has closed form expressions in frequency and space.
This provides an explicit characterization of $\smash{\psi_{j,k^{\nu},t^{\nu}}^{n-1}(\pmb{x})}$ and facilitates efficient numerical implementations, see Fig.~\ref{fig:conceptual}.

With Eq.~\ref{eq:psi_chain}, the projected signal along a direction $\nu$ is given by the \emph{local Fourier slice equation}
\begin{align}
  \label{eq:slice}
  f_{\nu}(\pmb{y})
  = \int_{\mathbb{R}_{\nu}} f(\pmb{x}) \, d\nu
  = \sum_{s \in \mathcal{I}} f_s \, \psi_{s}^{n-1,\nu}(\pmb{y})
\end{align}
that yields $f_{\nu}(\pmb{y})$ as wavelet reconstruction with the $(n-1)$-dimensional polar wavelets $\psi_{s}^{n-1,\nu}(\pmb{y})$.
The inverse Fourier transform in Eq.~\ref{eq:fourier_slice} is in our formulation thus replaced by a sum over wavelet coefficients, which can be implemented without the need for a further discretization, as usually required for the classical slice theorem.

With Eq.~\ref{eq:slice}, sparsity in the wavelet representation of a signal is readily exploited by restricting the sum to nonzero coefficients.
Furthermore, since Eq.~\ref{eq:psi_chain} preserves directionality
\setlength{\abovecaptionskip}{2pt}
\setlength{\columnsep}{4pt}
\setlength{\abovecaptionskip}{10pt}
and angular localization, depending on $\nu$, sparsity can  also be conserved.
By using that the $\psi_{s}^{n-1,\nu}(\pmb{y})$ are again wavelets, the projected reconstruction in Eq.~\ref{eq:slice} can also be performed locally, in space to obtain the projected signal over a sub-domain, or in frequency, to obtain a filtered version.
The computational complexity is thereby given by $\mathcal{O}(\vert \bar{\pmb{x}}_{\nu} \vert \, k \, \omega)$ and depends linearly on the size $\vert \bar{\pmb{x}}_{\nu} \vert$ of the region $\bar{\pmb{x}}_{\nu}$ onto which the signal is projected, the number $k$ of coefficients in the sparse signal representation, and the fraction $\omega$ of wavelets aligned with the projection direction.
The costs hence scale in a direct and intuitive manner with the size of the region of interest and the complexity of the signal with respect to the projection direction.

We numerically validate our local Fourier slice equation, demonstrating that it enables a local and sparse projection and with computational costs in correspondence with the theory.
The relevance of our Fourier slice equation for applications in optics is exemplified using tomographic reconstruction.
We show in particular how sparsity can be used as a ``magnifying lens'' to reconstruct with a higher resolution around a region of interest, thereby saving orders of magnitude in computation time and memory.

For concreteness we will restrict the following discussion to two- and three dimensions.
The conventions used in our work as well as some derivations are relegated to the appendix.

\subsection{Related Work}

The classical Fourier slice theorem~\cite{Bracewell1956} and the closely related Radon transform~\cite{Radon1917} have been used in a wide range of applications.
However, only a limited number of works combined them with wavelets or wavelet-like constructions.
Cand{\`e}s and Donoho~\cite{Candes2000} used curvelets as a sparsity prior to improve tomographic reconstruction from noisy measurements.
Later, Frikel~\cite{Frikel2013} improved upon their results and in particular considered limited angle tomography.
Garduno and co-workers~\cite{Garduno2011,Garduno2017} studied tomographic reconstruction using Haar wavelets, demonstrating that no improvement over classical approaches can be obtained using the algorithm they employed.
Shearlets, which are essentially a stereographic projection of polar wavelets, have been employed for sparse tomographic reconstruction using optimization~\cite{Garduno2017,Vandeghinste2013}, since sparsity is difficult to incorporate into algebraic approaches.
De Hoop et al.~\cite{DeHoop2009a} used a curvelet-like frame for tomographic reconstruction problems in geoscience, exploiting information about the relevant partial differential equation which we do not assume in our work.
Using concepts from compressed sensing, J{\o}rgensen et al.~\cite{Jorgensen2015a,Jorgensen2015} investigated how many measurements are required for optimization-based, sparse reconstructions.
We also exploit sparsity but using linear least squares reconstruction.
To our knowledge, the intrinsic connection between polar wavelets (including curvelets) and the geometry of the Fourier slice theorem, and that this results in a closed sequence of wavelets, has not been observed in the literature before.

Wavelet-like constructions defined in polar coordinates in the Fourier domain have been proposed in various forms over the years, e.g.~\cite{Simoncelli1995,Portilla2000,Candes1999b,Candes2005a}.
We build on the systematic framework recently proposed by Unser and co-workers~\cite{Unser2010,Unser2011,Ward2014}, which we refer to as polar wavelets~\cite{Lessig2018a}.

\section{A Local Fourier Slice Equation}
\label{sec:construction}

In this section we derive the local Fourier slice equation in Eq.~\ref{eq:slice}.
We will begin by briefly re-calling the construction of polar wavelets, which provides the basis for our work.
Then the two-dimensional case will be discussed before turning to the three-dimensional setting.

\subsection{Polar Wavelets}
\label{sec:construction:polar_wavelets}

\begin{figure}
  \centering
  \includegraphics[width=0.42\columnwidth]{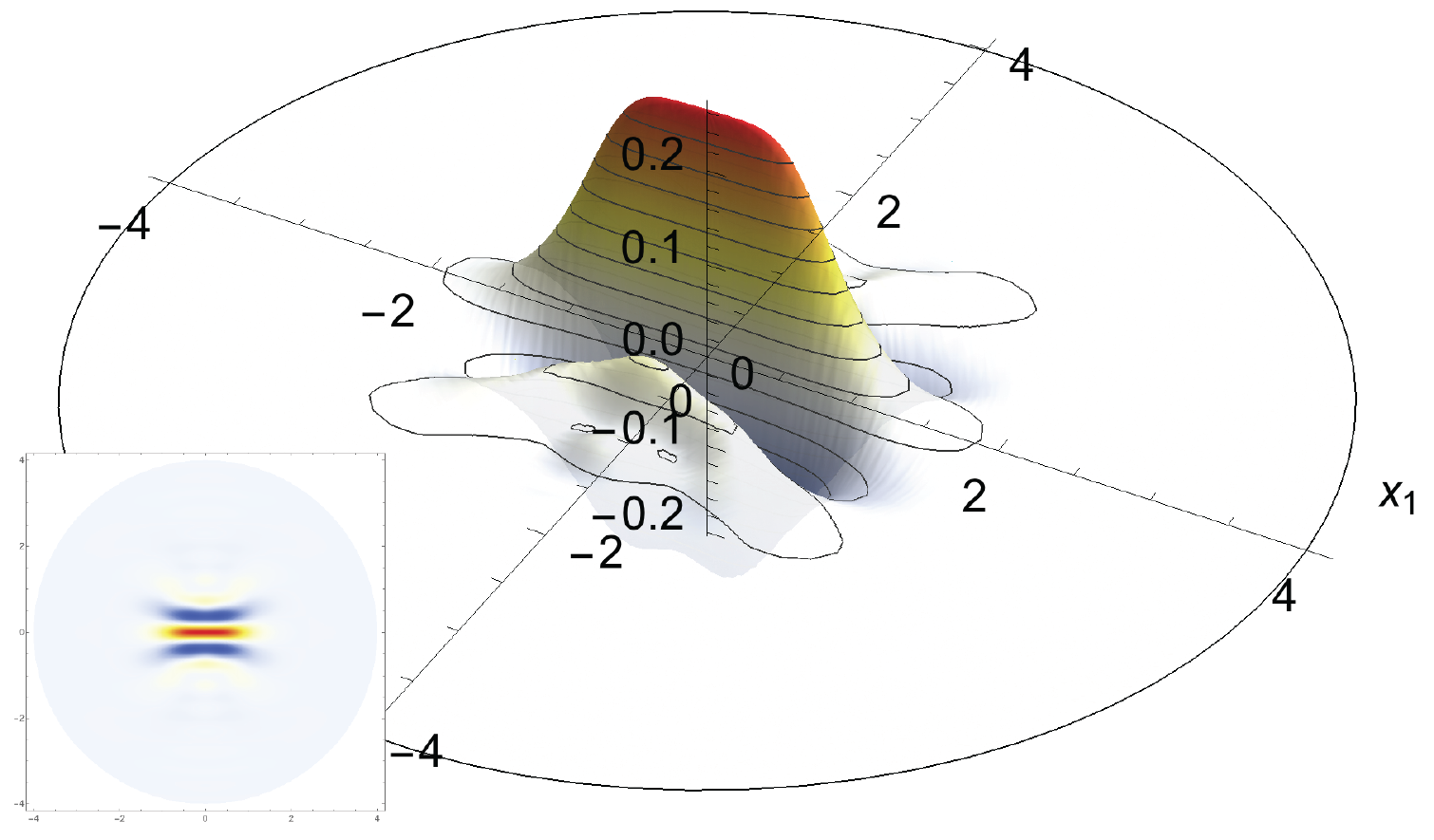}
  \hspace{0.2cm}
  \includegraphics[width=0.42\columnwidth]{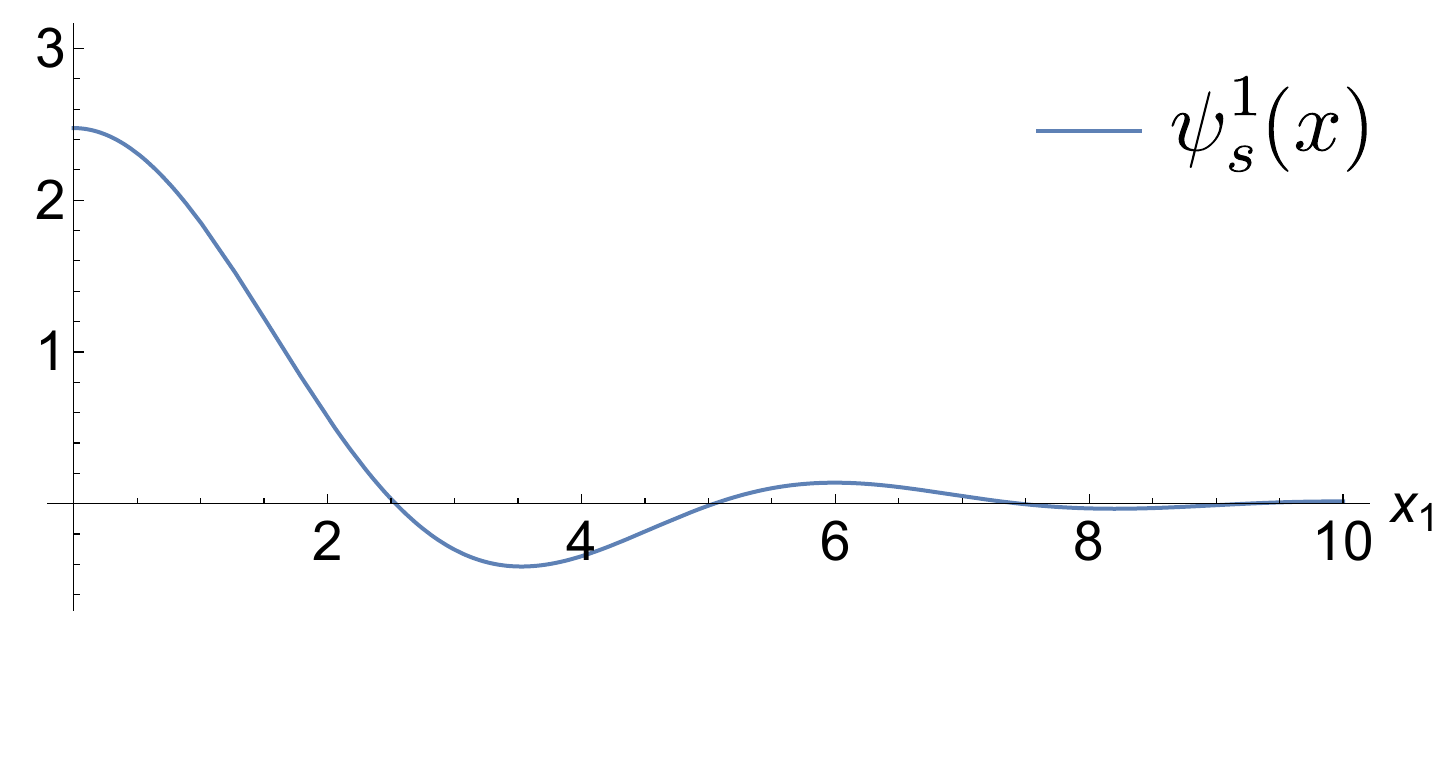}
  \caption{Directional polar wavelet $\psi_s^2(\pmb{x})$ in two dimensions (left) and its ``sliced'' counter-part $\psi_{0,0,0}^{1}(x_1) = h^1(\vert x_1 \vert)$ obtained by projecting along the $x_2$-axis (right).}
  \label{fig:zeta}
\end{figure}

Polar wavelets are defined in polar or spherical coordinates in the Fourier domain using a compactly supported radial window $\hat{h}( \vert \pmb{\xi} \vert)$, which controls the overall frequency localization, and an angular one $\hat{\gamma}(\theta_{\pmb{\xi}})$, which controls the directionality.
The mother wavelet is thus given by $\hat{\psi}(\pmb{\xi}) = \hat{\gamma}(\theta_{\pmb{\xi}}) \, \hat{h}( \vert \pmb{\xi} \vert )$ with the whole family of functions being generated by dilation, translation and rotation.

In two dimensions, the angular window is best described using a Fourier series.
A polar wavelet takes there hence the form
\begin{subequations}
  \label{eq:polarlets:2d}
\begin{align}
  \label{eq:polarlets:2d:hat}
  \hat{\psi}_s(\pmb{\xi})
  \equiv \hat{\psi}_{jkt}(\pmb{\xi})
  = \Big( \sum_{n} \beta_{j,n}^t \, e^{i n \theta_{\pmb{\xi}}} \Big) \, \hat{h}(2^{-j} \vert \pmb{\xi} \vert ) \, e^{-i \langle \pmb{\xi} , 2^j \pmb{k}\rangle}
\end{align}
with the $\beta_{j,n}^t$ controlling the angular localization.
In the simplest case $\beta_n = \delta_{n 0}$ and one has isotropic, bump-like wavelet functions.
In the spatial domain, the wavelets are given by
\begin{align}
  \label{eq:polarlets:2d:space}
  \psi_{s}(\pmb{x})
  \equiv \psi_{jkt}(\pmb{x})
  = \frac{2^j}{2\pi} \sum_{n} i^n \, \beta_{j,n}^t \, e^{i n \theta_{\pmb{x}}} \, h_n(2^{j} \vert \pmb{x} -\pmb{k}\vert )
\end{align}
  \end{subequations}
where $h_n(\vert \pmb{x} \vert)$ is the Hankel transform of $\hat{h}(\vert \pmb{\xi} \vert)$ of order $n$.
For $\hat{h}(\vert \pmb{\xi} \vert)$ we will employ the window proposed for the steerable pyramid~\cite{Portilla2000}, since $h_n(\vert \pmb{x} \vert)$ then has a closed form expression~\cite{Lessig2018a}.
When the wavelets in Eq.~\ref{eq:polarlets:2d}  are suitably augmented using scaling functions $\phi_{j,k}(\pmb{x})$ to represent a signal's low frequency part, with $\psi_{-1,k}(\pmb{x}) \equiv \phi_{0,k}(\pmb{x})$, the polar wavelets in Eq.~\ref{eq:polarlets:2d} provide a tight frame for $L_2(\mathbb{R}^2)$.
Hence any signal $f(\pmb{x}) \in L_2(\mathbb{R}^2)$ can be represented as
\begin{align}
  \label{eq:polarlets:2d:frame}
  f(\pmb{x})
  = \sum_{s \in \mathcal{I}} \Big\langle f(\pmb{y}) \, , \, \psi_s(\pmb{y}) \Big\rangle \, \psi_s(\pmb{x})
  = \sum_{j=-1}^{\infty} \sum_{k \in \mathbb{Z}^2} \sum_{t=1}^{N_j} \left\langle f(\pmb{y}) \, , \, \psi_{jkt}(\pmb{y}) \right\rangle \, \psi_{jkt}(\pmb{x})
\end{align}
and, although redundant, the frame affords most of the conveniences of an orthonormal basis.

Analogous to Eq.~\ref{eq:polarlets:2d:hat}, in three dimensions polar wavelets are defined by
\begin{align}
  \label{eq:polarlets:3d:hat}
  \hat{\psi}_{j,k,t}(\pmb{\xi})
  =  \hat{\gamma}_{j,t}\big(\bar{\pmb{\xi}} \big) \, \hat{h}(2^{-j} \vert \pmb{\xi} \vert) \, e^{-i \langle \pmb{\xi} ,2^j \pmb{k}\rangle}
  =  \sum_{l,m} \kappa_{lm}^{jt} \, y_{lm}(\bar{\pmb{\xi}}) \, \hat{h}(2^{-j} \vert \pmb{\xi} \vert) \, e^{-i \langle \pmb{\xi} , 2^j \pmb{k}\rangle}
\end{align}
where $\bar{\pmb{\xi}} = \pmb{\xi} / \vert \pmb{\xi} \vert$, the $y_{lm}(\bar{\pmb{\xi}})$ are spherical harmonics, and the coefficients $\kappa_{lm}^{jt}$ control the angular localization.
The wavelets in Eq.~\ref{eq:polarlets:3d:hat} have again closed form expressions in the spatial domain and they generate a tight frame for $L_2(\mathbb{R}^3)$, so that the analogue of Eq.~\ref{eq:polarlets:2d:frame} holds for all $f(\pmb{x}) \in L_2(\mathbb{R}^3)$.
We refer to the original works~\cite{Unser2013,Ward2014} and~\cite{Lessig2018a} for a more detailed discussion of polar wavelets.

\begin{figure}[t]
  \centering
  \hspace{0.05in}
  \includegraphics[trim={6.5cm, 10.5cm, 4cm, 6cm},clip, width=0.36\textwidth]{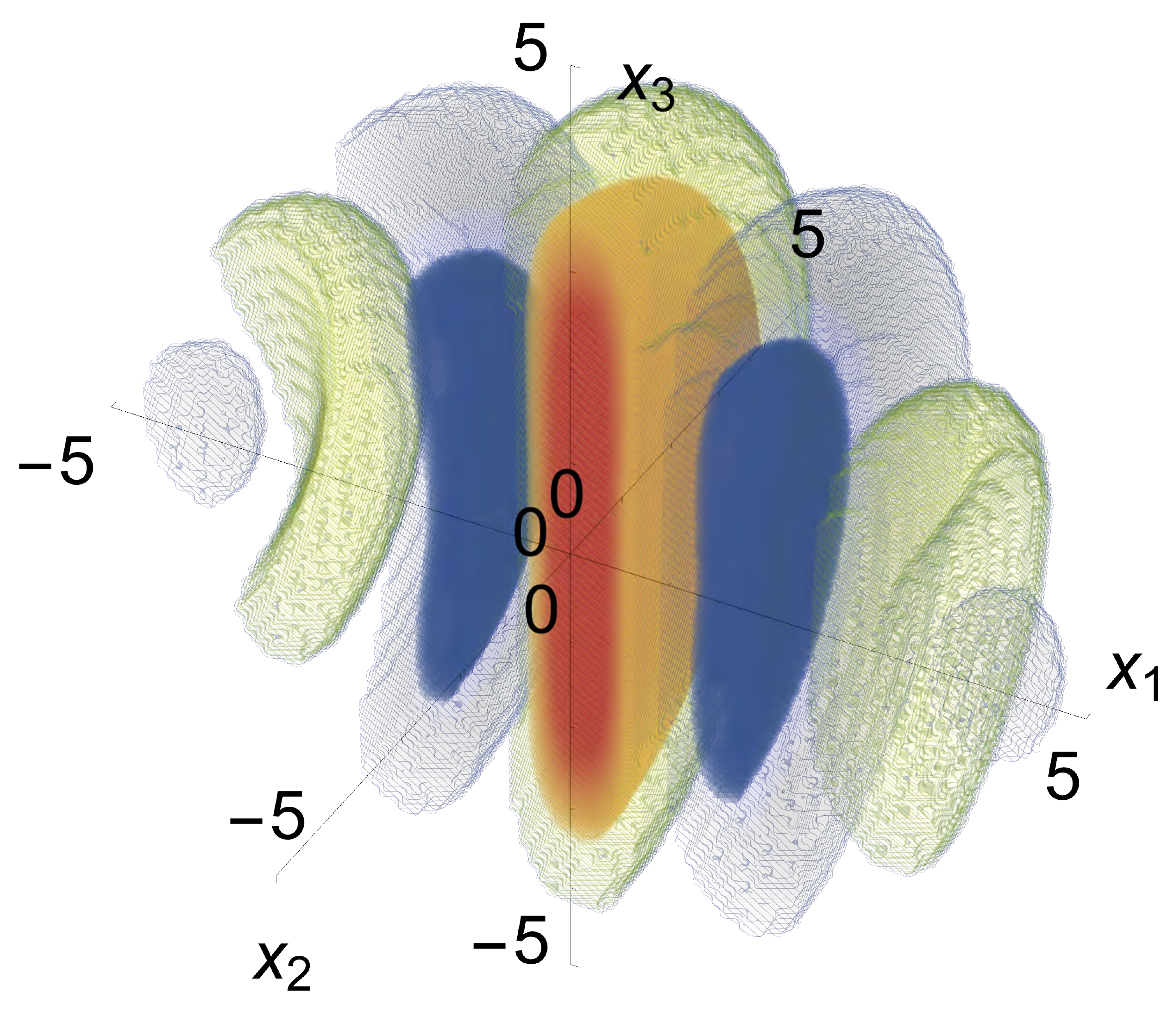}
  \hspace{0.30in}
  \includegraphics[width=0.40\columnwidth]{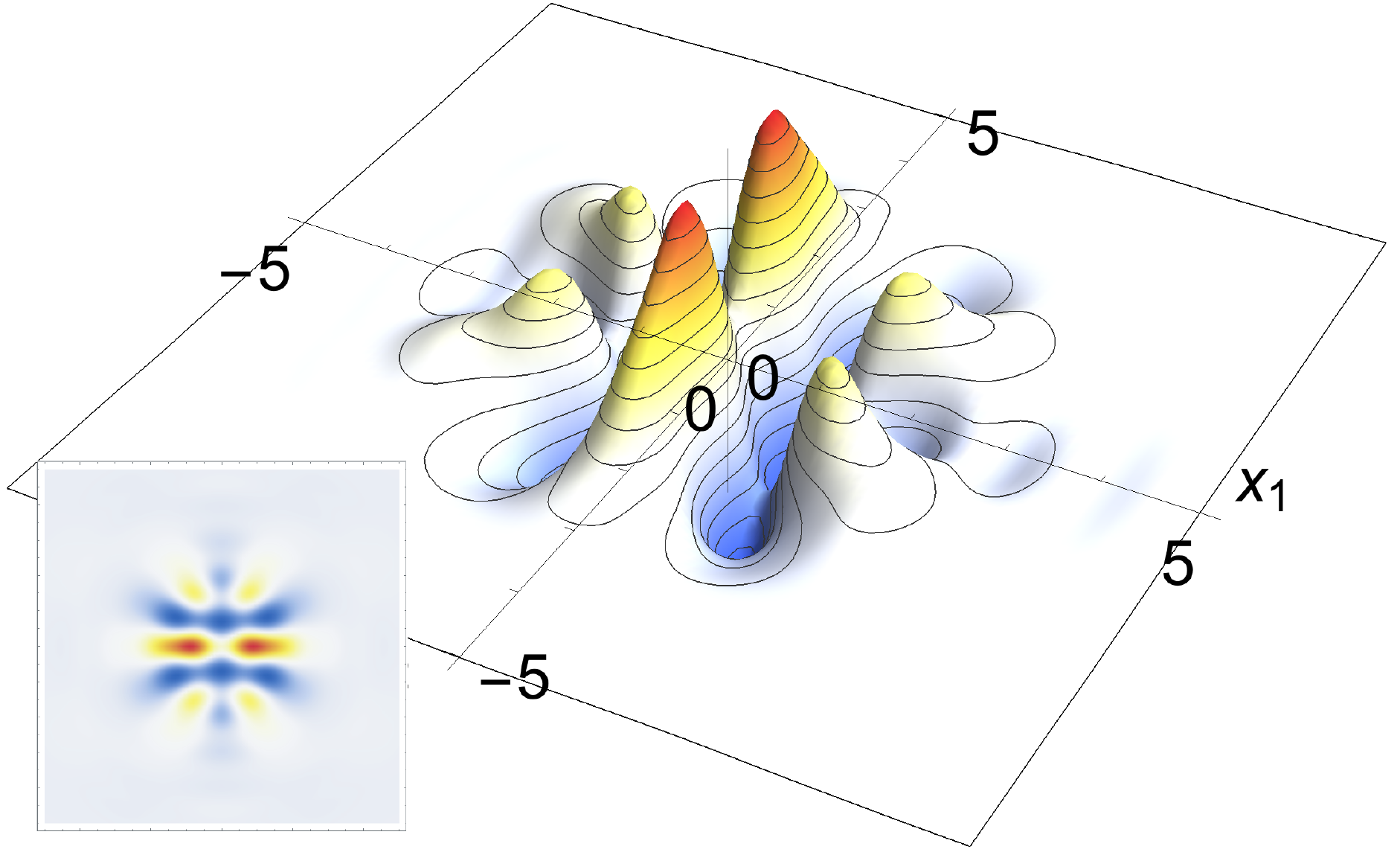}
  \caption{Directional polar wavelet $\psi_s^3(\pmb{x})$ in $\mathbb{R}^3$ and its projection $\psi_s^2(x_{12})$ along the $x_3$ axis, which is a two-dimensional polar wavelet. Note how the orientation of $\psi_s^3(\pmb{x})$ is essentially preserved under projection. This is critical for the conservation of sparsity.
  }
  \label{fig:zeta:2d}
\end{figure}

\subsection{Local Fourier Slice Equation in the Plane}
\label{sec:construction:slice:2d}

In the plane and when $\nu$ is in the directions of the $x_2$-axis, the classical Fourier slice theorem is easily established.
Writing $f(\pmb{x})$ as its inverse Fourier transform we have for the projection
\begin{subequations}
  \label{eq:slice:classical}
\begin{align}
  f_2(x_1)
  &= \frac{1}{2\pi} \int_{\mathbb{R}_{x_2}} \int_{\mathbb{R}_{\pmb{\xi}}^2} \hat{f}(\xi_1,\xi_2) \, e^{i \langle \pmb{\xi} , \pmb{x} \rangle} \, d\pmb{\xi} \, dx_2 .
\end{align}
Since the Fourier transform $\hat{f}(\xi_1,\xi_2)$ does not depend on $x_2$, the integral over $\mathbb{R}_{x_2}$ only involves $e^{i \langle \pmb{\xi} , \pmb{x} \rangle}$, yielding $e^{i \langle \xi_1 , x_1 \rangle} \, \delta(\xi_2)$. Thus
\begin{align}
  \label{eq:slice:classical:4}
  f_2(x_1)
  &= \frac{1}{2\pi} \int_{\mathbb{R}_{\xi_1}} \hat{f}(\xi_1,0) \, e^{i \langle \xi_1 , x_1 \rangle} \, d\xi_1
\end{align}
\end{subequations}
which is the Fourier slice theorem.
The general result, for an arbitrary axis of integration, follows by the covariance of the Fourier transform.

With $f(\pmb{x})$ in Eq.~\ref{eq:slice:classical:4} given in its polar wavelet representation,
\begin{subequations}
  \label{eq:fourier_slice:local:der}
\begin{align}
  \label{eq:fourier_slice:local:der:1}
  f_2(x_1)
  &= \frac{1}{2\pi} \int_{\mathbb{R}_{\xi_1}} \left( \sum_{s} f_s \, \hat{\psi}_s(\xi_1,0) \right) \, e^{i \langle \xi_1 , x_1 \rangle} \, d\xi_1 .
\end{align}
Using linearity and with the definition of the polar wavelets in Eq.~\ref{eq:polarlets:2d:hat} we obtain
\begin{align}
  \label{eq:fourier_slice:local:der:2}
  f_2(x_1)
  &= \frac{1}{2\pi} \sum_{s} f_s \,
  \hat{\gamma}_s(0) \! \int_{\mathbb{R}_{\xi_1}} \! \! \! \hat{h}(2^{-j_s} \vert \pmb{\xi} \vert ) \, e^{-i \langle \xi_1 , k_1^s \rangle}  \, e^{i \langle \xi_1 , x_1 \rangle} \, d\xi_1
\end{align}
\end{subequations}
where, through the ``slicing'', the angular window $\hat{\gamma}(\theta_{\pmb{\xi}})$ no longer depends on the integration variable and only needs to be evaluated at $\theta_{\pmb{\xi}} = 0$.
The remaining integral in Eq.~\ref{eq:fourier_slice:local:der:2} is a one-dimensional Fourier transform with translation factor $e^{-i \langle \xi_1 , k_1^s \rangle}$.
By defining
\begin{align}
  \label{eq:zeta:2d:aligned}
  \psi_{s}^1(x_1) = \frac{1}{2\pi} \, \hat{\gamma}_s(0)
   \underbrace{\int_{\mathbb{R}_{\xi_1}} \! \! \! \hat{h}\big(2^{-j_s} \vert \pmb{\xi} \vert \big) \, e^{-i \langle \xi_1 , k_1^s \rangle}  \, e^{i \langle \xi_1 , x_1 \rangle} \, d\xi_1}_{\displaystyle h^1( 2^j x_1 - k_1^s) }
\end{align}
we recover the local Fourier slice equation in Eq.~\ref{eq:slice}.
For the radial window $\hat{h}(\vert \pmb{\xi} \vert)$ of the steerable pyramid~\cite{Portilla2000}, the profile $h^1(\vert x \vert)$ in Eq.~\ref{eq:zeta:2d:aligned} has, in fact, a closed form expression,
\begin{align}
  \label{eq:fourier_slice:local:2}
   h^1(x) &= - \frac{i \pi}{8}  \!
   \big(E_{c}(-z) + E_{c}(z) + 4 E_c(4 z) + 4 E_c(-z)\big)
\end{align}
where $E_{c}(z) = E_{c_+}(z) - E_{c_-}(z)$ and $E_n(\cdot)$ is the exponential integral function, $z = i \pi x / 4$, and $c_{\pm} = \pm (i \pi) / \mathrm{log}(4)$, see Fig.~\ref{fig:zeta} for a plot.

Polar wavelets are covariant under rigid body motions~\cite{Azencott2009}.
The above derivation hence immediately implies the result for an arbitrary projection direction $\nu$ since we can first rotate the original polar wavelet representation so that $\nu$ is aligned with $x_2$, which amounts to rotating the grid over which the wavelets are defined, then applying the above result, and finally rotating the projected signal back onto $x_v$, the axis orthogonal to $\nu$.
We thus have the following result.

\begin{figure}
  \centering
  \includegraphics[width=1.0\textwidth]{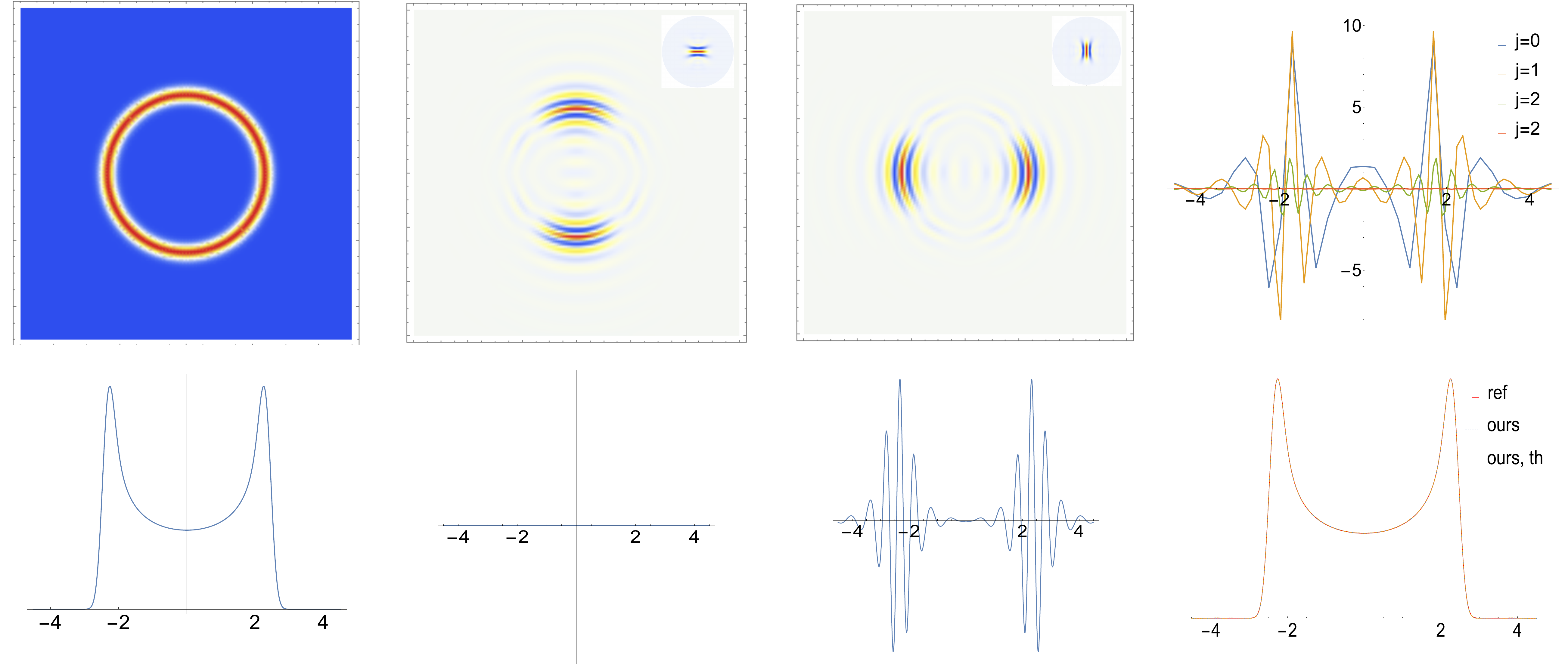}
  \caption{\emph{Left:} Original signal and its projections onto the $x_1$-axis (horizontally aligned). \emph{Center columns:} Basis function coefficients (top) and resulting contribution to the projected signal (bottom) for horizontal and vertical orientation of the wavelet (insets) \emph{Right, top:} Coefficients $\smash{f_{j,k_1,t}^1}$ of the projected signal. They are small away from the signal, providing an example for the conservation of sparsity under ``slicing''. \emph{Right, bottom:} projected signal obtained using the local Fourier slice equation with all (dotted, ``ours'') and only the $5\%$ largest coefficients (dashed, ``ours, th'').}
  \label{fig:rim}
\end{figure}

\begin{proposition}
  \label{prop:2d}
  Let $f(x) \in L_1(\mathbb{R}^2) \cap L_2(\mathbb{R}^2)$ and $\{ \psi_s(x) \}_{s \in \mathcal{I}}$ be a Parseval tight polar wavelet frame for $L_2(\mathbb{R}^2)$.
  Then the projection of $f(x)$ along direction $\nu \in S^1$, with polar coordinate $\theta_{\nu}$, is given by the local Fourier slice equation in Eq.~\ref{eq:slice} with
  \begin{align}
    \label{eq:prop:2d}
    \psi_s^{1,\nu}(x_{\nu}) = \psi_{j_s,k_s^{\nu}}^{1,\nu}(x_{\nu}^1) = \frac{1}{2\pi} \, \hat{\gamma}(\theta_\nu) \, h^1( 2^j x_{\nu} - k_s^{\nu} )
  \end{align}
  and $k_s^{\nu}$ being the projection of $k_s$ onto the line orthonormal to $\nu$.
\end{proposition}

Except when $\nu$ is along an axis, the wavelets $\psi_s^{1,\nu}(x_{\nu})$ are no longer equi-spaced but positioned at the irregular locations $k_s^{\nu}$.
Eq.~\ref{eq:slice} nonetheless holds by construction.

\begin{remark}
A derivation analogous to those in Eq.~\ref{eq:fourier_slice:local:der} can also be performed for classical tensor product wavelets.
However, for every direction $\nu$ one then has a different projected wavelet that is spread across multiple scales $j$ and that does not have a closed form expression or simple description.
How an efficient implementation would be possible is hence unclear.
Also, directional sparsity could not be exploited, since the wavelets are not directionally localized.
The latter one would be possible with contourlets~\cite{Do2005a} and shearlets~\cite{Labate2005} but with these one would only approximately obtain a $1$-dimensional wavelet and, to our knowledge, no closed form expression for it would be available.
\end{remark}

\subsection{Local Fourier Slice Equation in Space}
\label{sec:construction:slice:3d}

In three dimensions, two different projections are possible. Projecting along one axis only, which corresponds to the X-ray transform, and along two axes, so that one again obtains a one-dimensional signal.
The derivations proceed in both cases analogously to the two-dimensional setting we discussed in detail in the previous sub-section and they thus have been relegated to Appendix~\ref{sec:appendix:3d}.
We summarize the results in the following propositions.

\begin{proposition}
  \label{prop:3d:1axis}
   Let $f(\pmb{x}) \in L_1(\mathbb{R}^3) \cap L_2(\mathbb{R}^3)$ and $\{ \psi_s(\pmb{x}) \}_{s \in \mathcal{I}}$ be a Parseval tight polar wavelet frame for $L_2(\mathbb{R}^3)$, as defined in Eq.~\ref{eq:polarlets:3d:hat}.
  Then the projection of $f(\pmb{x})$ along direction $\nu \in S^2$ is given by the local Fourier slice equation in Eq.~\ref{eq:slice} with
  \begin{subequations}
    \label{eq:prop:3d:1axis}
  \begin{align}
    \label{eq:prop:3d:1axis:1}
     \psi_s^{2,\nu}( \pmb{x}_{\nu} )
      =  \sum_{m} i^m \, \beta_{m}^{jt,\nu} \, e^{i \, m \, \phi_{\pmb{x}_{\nu}}} \, h_m\Big( 2^{-j} \vert \pmb{x}_{\nu} - \pmb{k}_s^{\nu} \vert \Big)
    \end{align}
    and angular localization coefficients
    \begin{align}
      \label{eq:prop:3d:1axis:2}
      \beta_{m}^{jt,\nu} = \sum_l C_{lm} \, \left( \sum_{m'} W_{lm}^{m'}(\nu) \kappa_{lm'}^{j t}  \right) \, P_{lm}\big(\pi / 2\big)
    \end{align}
    \end{subequations}
    where the $W_{lm}^{m'}(\nu)$ are Wigner-D matrices, aligning $\nu$ with the $\xi_3$ axis, $\pmb{k}_s^{\nu}$ is the projection of $\pmb{k}_s$ onto the plane with normal $\nu$, and $h_m(\cdot)$ is the inverse Hankel transform of $\hat{h}(\cdot)$.
\end{proposition}

\begin{figure}
  \centering
  \includegraphics[width=0.49\columnwidth]{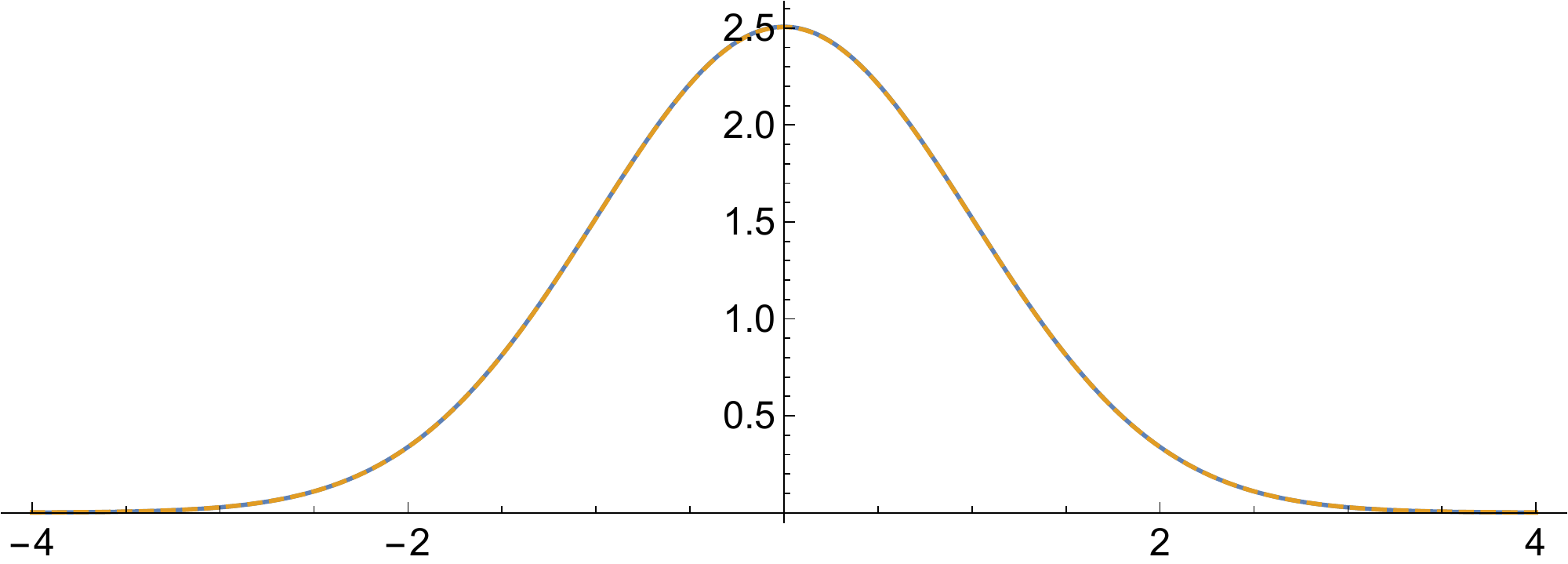}
  \includegraphics[width=0.49\columnwidth]{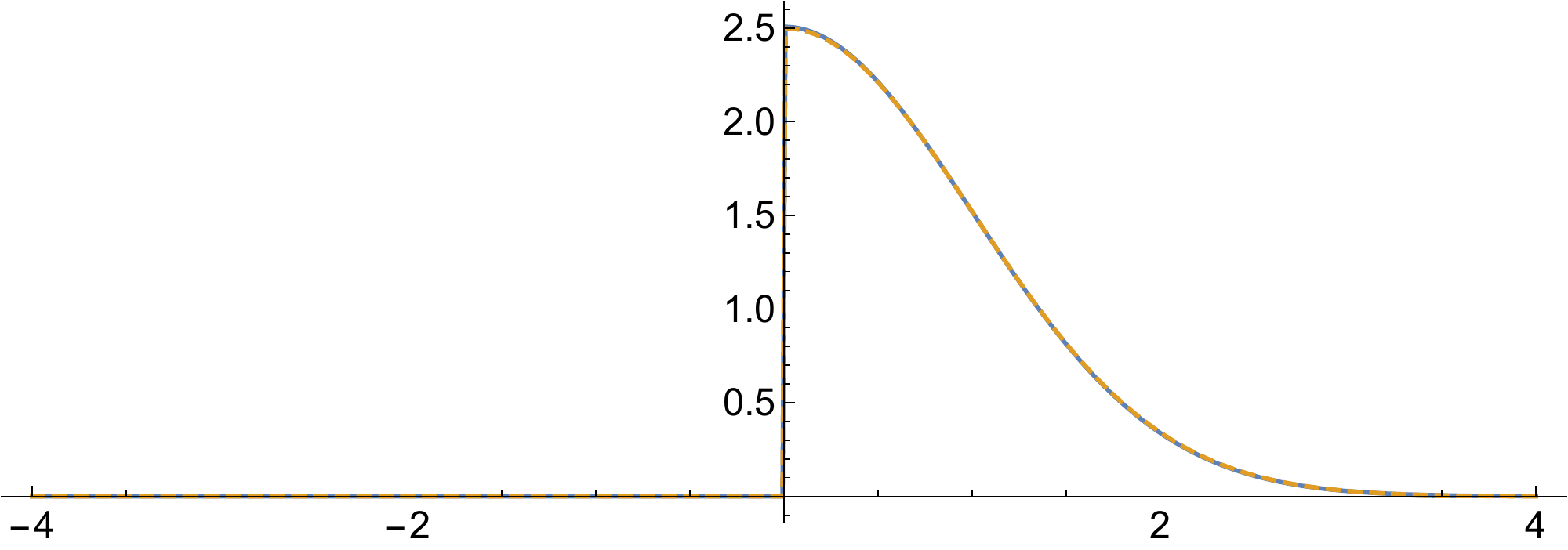}
  \caption{Local Fourier slice equation applied to a $2D$ Gaussian over the entire $x_1$ axis and localized only over the  positive one (blue, reference solution; yellow, local Fourier slice equation). The local reconstruction required $55\%$ of the time for the one over the entire axis.}
  \label{fig:gaussian}
\end{figure}

Comparing Eq.~\ref{eq:prop:3d:1axis} to Eq.~\ref{eq:polarlets:2d} we see that the one-dimensional projection of a three-dimensional polar wavelet yields a two-dimensional one and that the angular localization is preserved, to the extent possible, since the $\beta_{m}^{jt,\nu}$ are obtained from the original angular coefficients $\kappa_{lm'}^{j t}$. In particular, the projection of an isotropic polar wavelet in $\mathbb{R}^3$ yields an isotropic one in $\mathbb{R}^2$.
Proposition~\ref{prop:3d:1axis} provides the frequency representation of the projected wavelet.
But since it is a two dimensional polar wavelet, the spatial representation is immediately given by Eq.~\ref{eq:polarlets:2d:space}.
An example is shown in Fig.~\ref{fig:zeta:2d}.

For the projection along two axes we have the following result.

\begin{proposition}
  \label{prop:3d:2axis}
   Let $f(\pmb{x}) \in L_1(\mathbb{R}^3) \cap L_2(\mathbb{R}^3)$ and $\{ \psi_s(\pmb{x}) \}_{s \in \mathcal{I}}$ be a Parseval tight polar wavelet frame for $L_2(\mathbb{R}^3)$ as defined in Eq.~\ref{eq:polarlets:3d:hat}.
  Then the projection of $f(\pmb{x})$ onto the axis $\pmb{x}_{\nu}$ is given by the local Fourier slice equation in Eq.~\ref{eq:slice} with
  \begin{align}
    \label{eq:prop:3d:2axis}
     \psi_s^{1,\nu}( x_{\nu} )
      &= \hat{\gamma}_{s}\big( R_{\nu} \xi_{3} \big) \, h^1( \vert x_{\nu} - k_s^{\nu} \vert)
  \end{align}
  and $k_s^{\nu}$ is the projection of $k_s$ onto the $x_{\nu}$ axis.
\end{proposition}

Comparing Eq.~\ref{eq:prop:3d:2axis} to Eq.~\ref{eq:prop:2d} we see that the projection onto one axis yields the same $1$-dimensional wavelet we obtained in Proposition~\ref{prop:2d} for the projection in $\mathbb{R}^2$.

\subsection{Conservation of Sparsity}
\label{sec:construction:sparsity}

To understand the effect of the local Fourier slice equation on sparsity, we begin with the projection along the $x_2$-axis in $\mathbb{R}^2$.
Since the projection is aligned with the translation grid, it becomes
\begin{align}
  \label{eq:2d:projection:aligned}
  f^1(x_1)
  = \sum_{s \in \mathcal{I}} f_s \, \psi_{s}^1(x_1)
  = \sum_{j,k_1,t} \Big( \sum_{k_2} f_{jkt} \Big)  \, \psi_{j,k_1,t}^1(x_1)
  = \sum_{j,k_1,t}  f_{j,k_1,t}^1 \, \psi_{j,k_1,t}^1(x_1)
\end{align}
with the sum over $k_2$ being decoupled and the $\psi_{s}^1(x_1)$ equispaced along the $x_1$ axis.
The wavelet coefficients of the projected signal, $f_{j,k_1,t}^1 = \sum_{k_2} f_{jkt}, $ inherit the space-frequency localization of the original representation since $\smash{f_{j,k_1,t}^1}$ is a superposition of the two-dimensional coefficients $f_{j,k,t}$ in the same frequency band $j$, for the same location $k_1$, and the same orientation $t$.
Hence, when the modulus of all $f_{j,k,t}$ is small then so is those of $f_{j,k_1,t}^1$.
However, sparsity can also be generated, when then sum over $k_2$ becomes small through cancellation, and it can be destroyed, when the $f_{jkt}$ accumulate to a non-negligible value.
We leave a systematic analysis of these cases to future work; existing results in this direction can be found in~\cite{Quinto1993,Quinto2007}.

\begin{figure}
  \centering
    \includegraphics[width=0.275\textwidth]{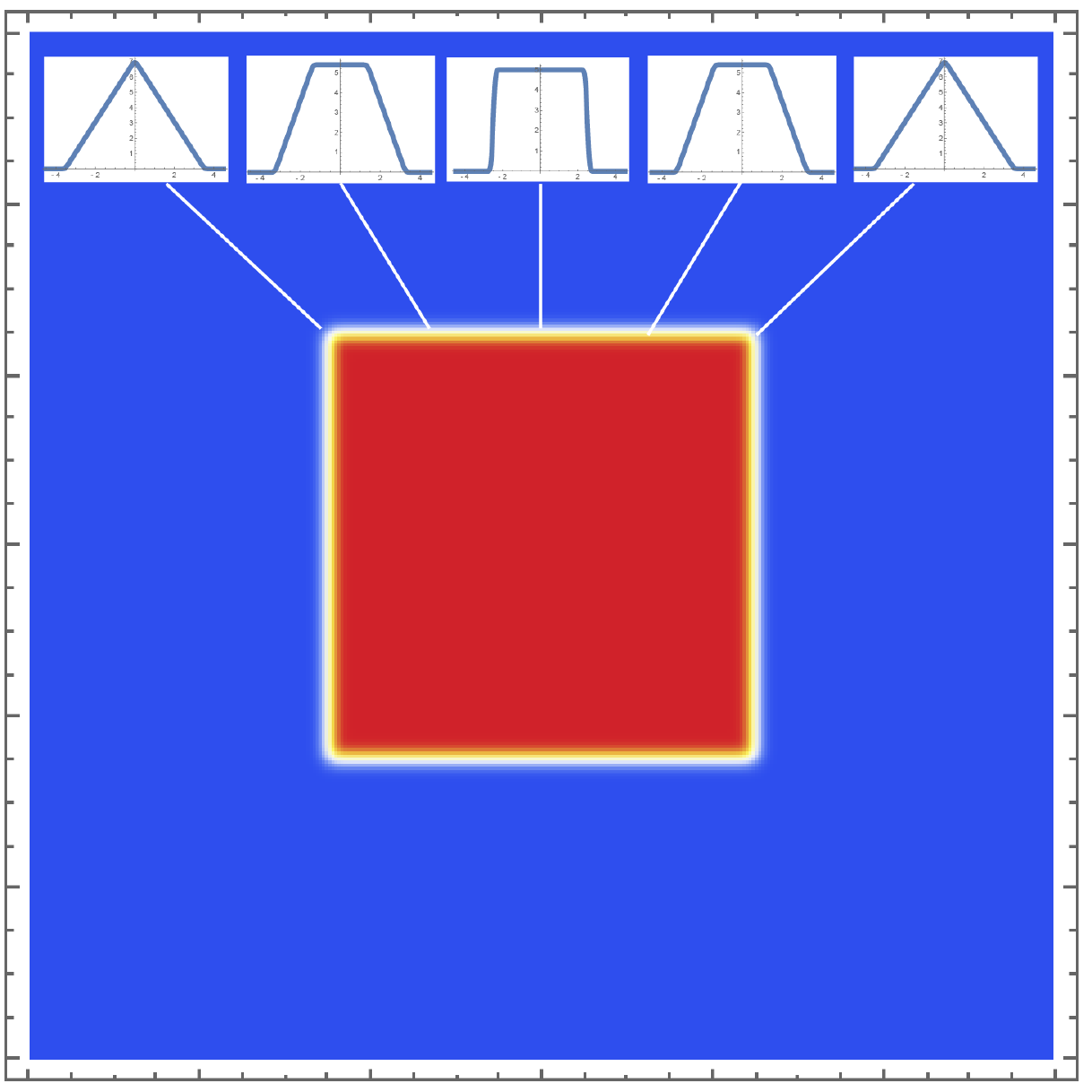}
    \hspace{0.2cm}
    \includegraphics[width=0.54\textwidth]{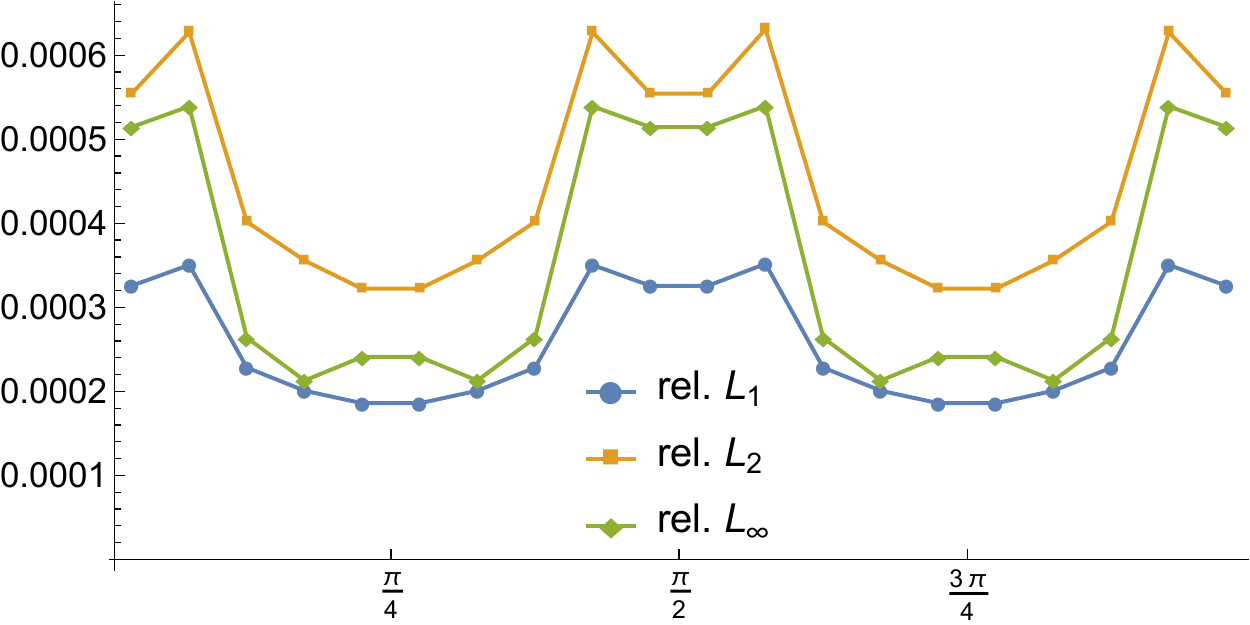}
  \caption{Relative $L_1$, $L_2$ and $L_{\infty}$ error (w.r.t. the respective norm of the original signal) of the Fourier slice equation for the box signal (left) as a function of the projection direction (projected signals as insets). Although the error fluctuates, it is overall independent of the direction. Variations result from changes in the regularity of the projection.}
  \label{fig:square}
\end{figure}

The contribution of $\smash{f_{j,k_1,t}^1}$ to the projected signal $f^1(x_1)$ will also be negligible when the value of the corresponding wavelet $\smash{\psi_{j,k_1,t}^1(x_1)}$ at the location $x_1$ is small.
By the spatial localization of the wavelet $\smash{\psi_{j,k_1,t}^1(x_1)}$ around $k_1$ this is true when $\vert k_1 - x_1 \vert \gg 2^{-j}$.
However, it also holds when the two-dimensional polar wavelet $\smash{\psi_{jkt}(x)}$ is not aligned with the projection direction.
Then the value of $\hat{\gamma}_s(0)$ will be small, which implies that the $x_2$-axis is along a direction where $\psi_s(\pmb{x})$ has vanishing moments.
An example is shown in Fig.~\ref{fig:rim} where $f(\pmb{x})$ is a thin annulus with radius $2.5$ centered at the origin.
Geometrically, the projection will have ``bumps'' around $x_1 = \pm 2.5$, since there one integrates along the rim, and it will be small around the origin, where one integrates orthogonal to it, cf. Fig.~\ref{fig:rim} left.
In a polar wavelet representation,  the coefficients $f_{jkt}$ will all have approximately the same nonzero magnitude when the wavelets are locally aligned with the annulus and be negligible for all other orientations.
In particular, around $x_1 = 0$ the coefficients $f_{jkt}$ will be significant only for horizontally oriented wavelets.
But then $\gamma_s(\theta_{x_2})$ vanishes and one obtains no contribution to the projection.
In contrast, around $x_1 = \pm 2.5$ the directional wavelets with a non-negligible coefficient are vertically oriented so that $\gamma_s(\theta_{x_2})$ is large and hence one obtains a significant contribution.
In Fig.~\ref{fig:rim}, right, one also sees the conservation of sparsity since for the region to the left and right of the signal the sparse representation of the two-dimensional signal is directly mapped to a sparse representation of the projected one.

The above discussion on the conservation of sparsity carries over to an arbitrary $\nu \in S^1$ by the covariance of polar wavelets and it applies with natural modifications to the situation in $\mathbb{R}^3$.
We leave a quantitative analysis to future work.

\subsection{Computational Complexity}
\label{sec:construction:costs}

In two dimensions, the computational costs of our local Fourier slice equation are given by $\mathcal{O}(\vert \bar{x}_{\nu} \vert \, k \, \omega)$ where $\vert \bar{x}_{\nu} \vert$ is the length of the projection region $\bar{x}_{\nu}$, since reconstruction of $f_{\nu}(y)$ typically amounts to determining it on a set of points whose total number is controlled by $\vert \bar{x}_{\nu} \vert$; $k$ is the number of coefficients in the sparse representation or the number of coefficients to be considered for projection, e.g. when not all levels are used, as follows from Eq.~\ref{eq:slice}; and $\omega$ is the fraction of basis functions aligned with the projection direction, which follows immediately from the factor $\hat{\gamma}(\theta_\nu)$ in $\psi_s^{1,\nu}(\bar{x}_{\nu})$ in Proposition~\ref{prop:2d}, see also Sec.~\ref{sec:construction:sparsity}.
The last two factors, which might depend on $j$, for example when the number of orientations is $j$-dependent, control the cardinality of the sum in the local Fourier slice equation in Eq.~\ref{eq:slice}.
The first factor determines how often the sum needs to be evaluated.
They all hence linearly affect the computational costs so that $\mathcal{O}(\vert \bar{x}_{\nu} \vert \,k\, \omega)$, with the constant being controlled by the cost of evaluating the projected wavelets $\psi_s^{n-1,\nu}(x_{\nu})$ and the density of the points over which $f_{\nu}(y)$ is determined.
In the next section we will see that this analysis indeed accurately describes the computational costs in our experiments.
$\mathcal{O}(\vert \bar{x}_{\nu} \vert \,k\, \omega)$ generalizes naturally to higher dimensions and $\vert \bar{x}_{\nu} \vert$ becomes then the area or volume of the region on which the projection is sought.

Our analysis assumes that the wavelet representation of the input signal is already available, e.g. because it is stored in a compressed form and different slices are to be determined as required.
If this is not the case than the fast wavelet transform can be employed to compute the wavelet representation in $O(n)$ time where $n$ is the number of samples in the input signal.

\section{Numerical Experiments}
\label{sec:experiments}

In this section we provide experimental results for our local Fourier slice equation.
We begin with basic experiments to provide insight into its fundamental behavior.
In Sec.~\ref{sec:experiments:tomography} we then discuss its application to tomographic reconstruction.
The code implementing the experiments is provided in the supplementary material~\cite{Lessig2018_slicecode} and implementation details are provided in Sec.~\ref{sec:appendix:implementation}.
We also provide information on absolute computation times but these should be considered as preliminary since we used a prototyping language.

\subsection{Basic experiments}
\label{sec:experiments:basic}

\begin{figure}
  \centering

  \centering
  \includegraphics[width=0.24\textwidth]{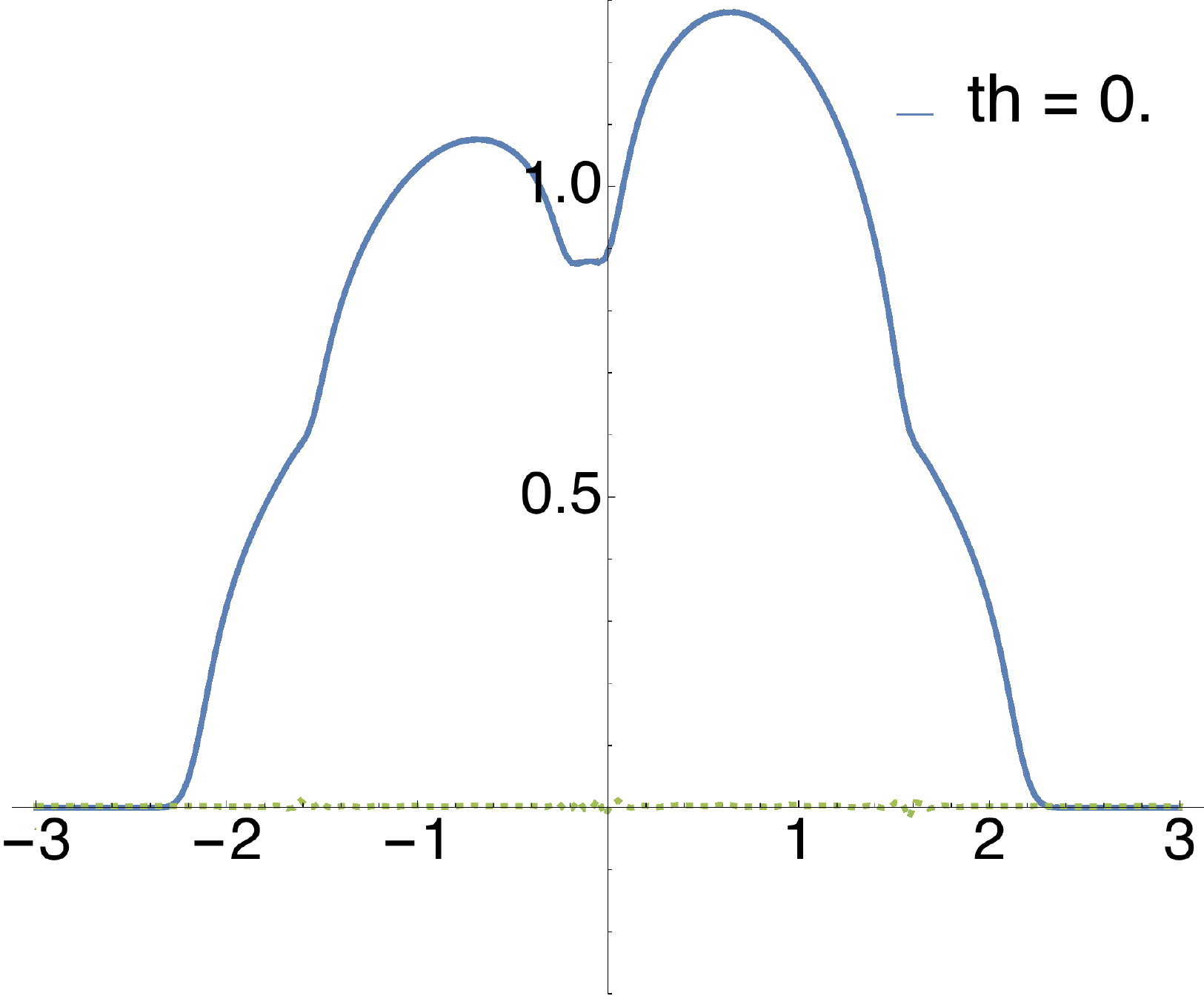}
  \includegraphics[width=0.24\textwidth]{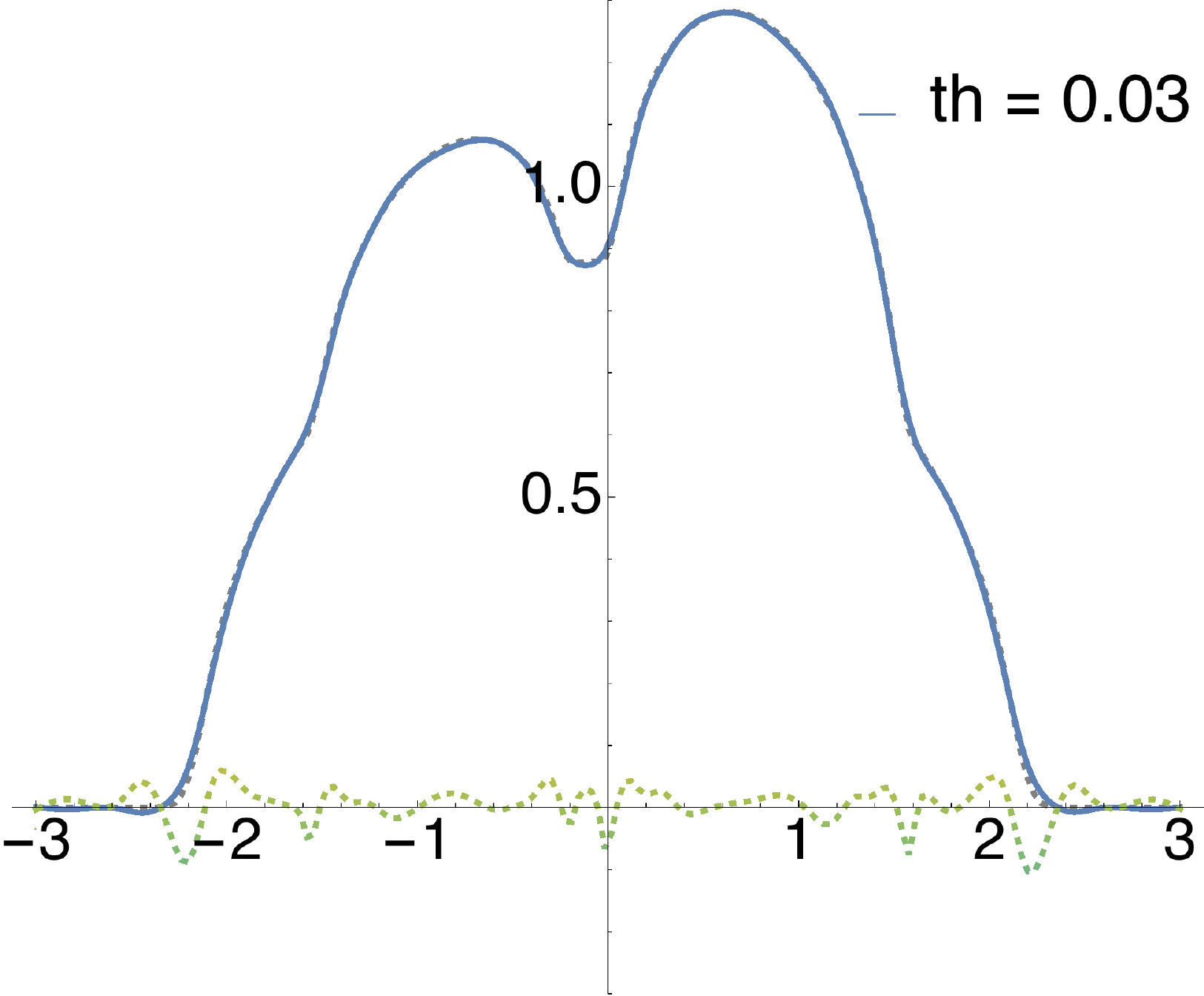}
  \includegraphics[width=0.24\textwidth]{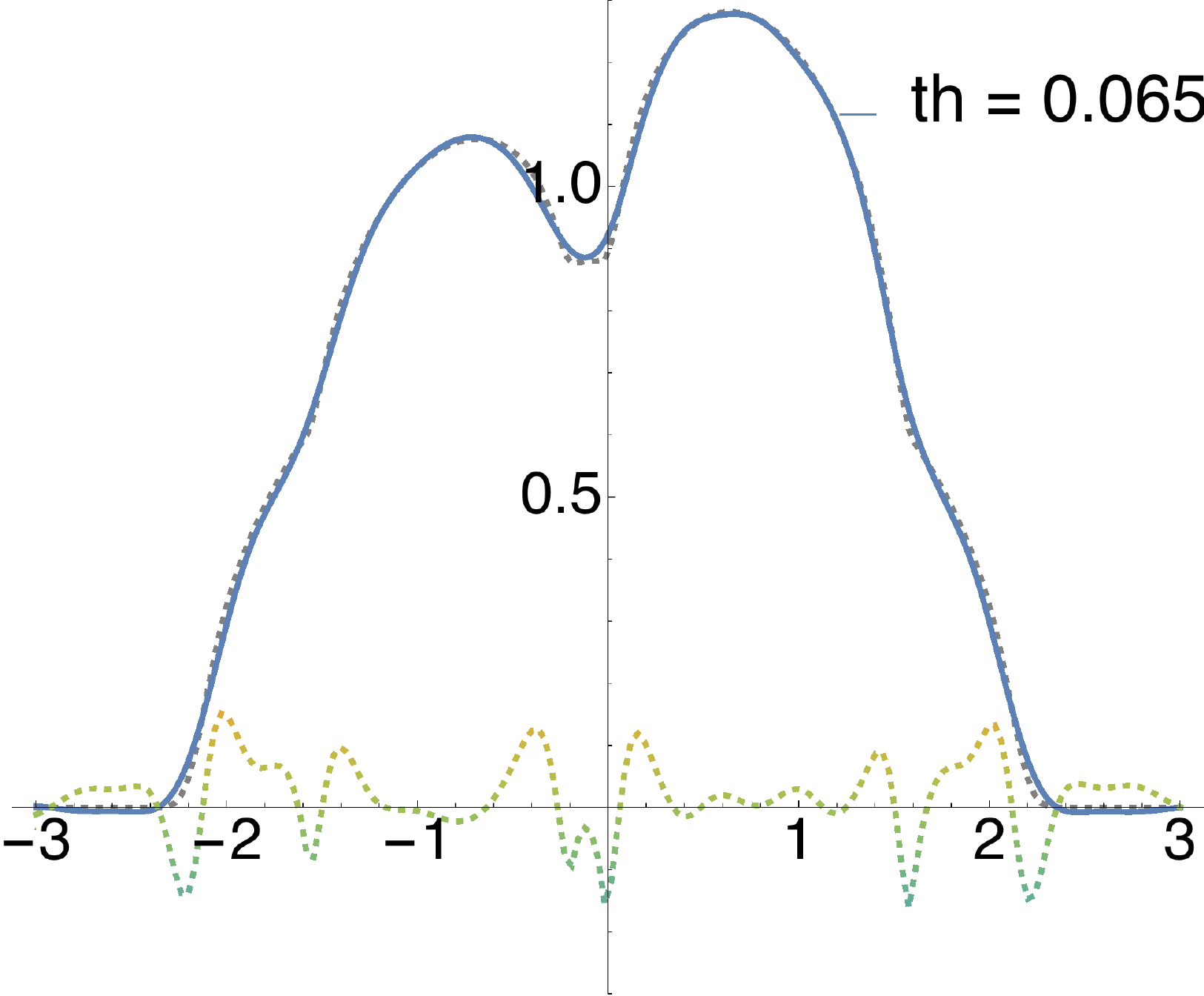}
  \includegraphics[width=0.24\textwidth]{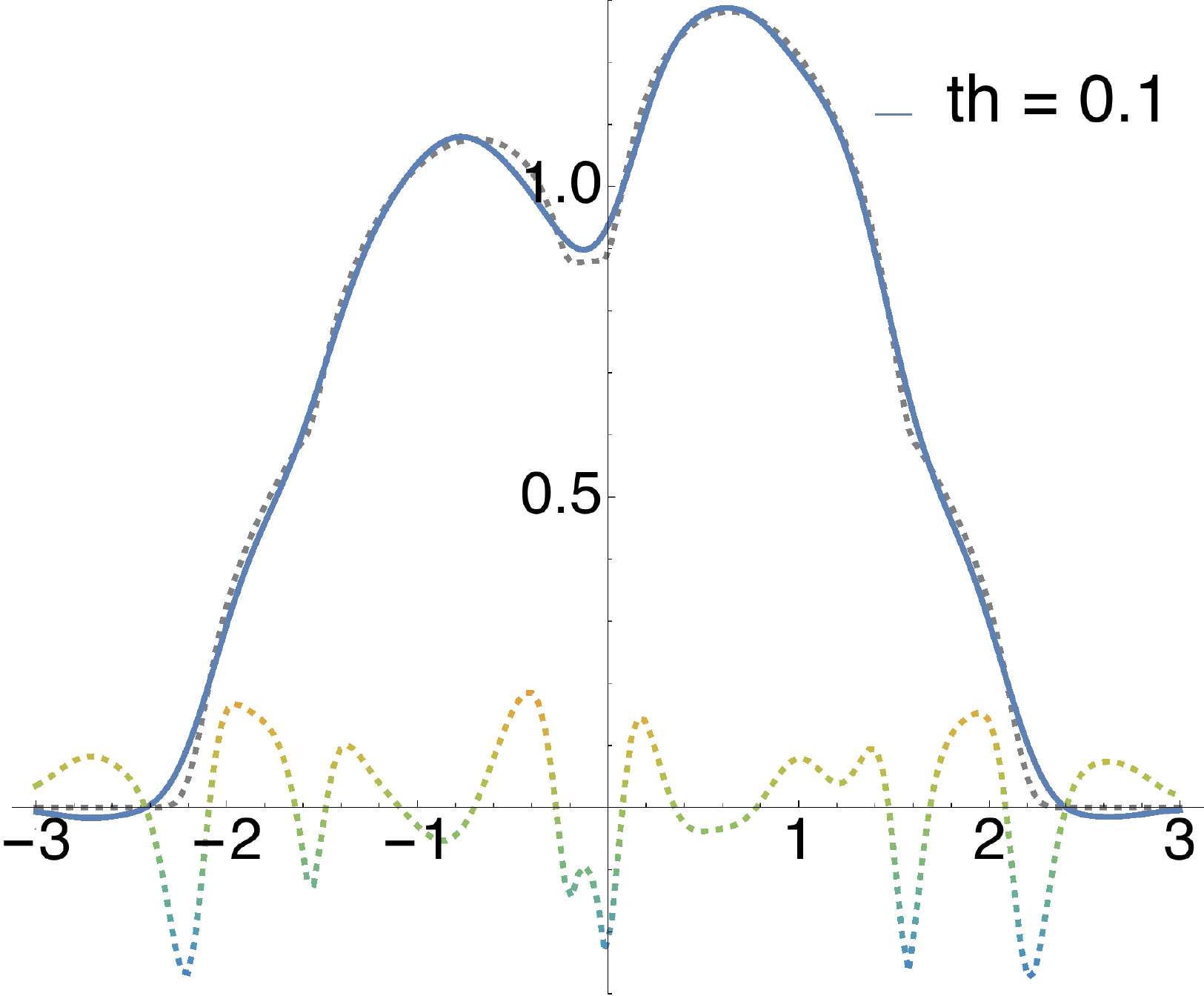}
  \includegraphics[width=0.475\textwidth]{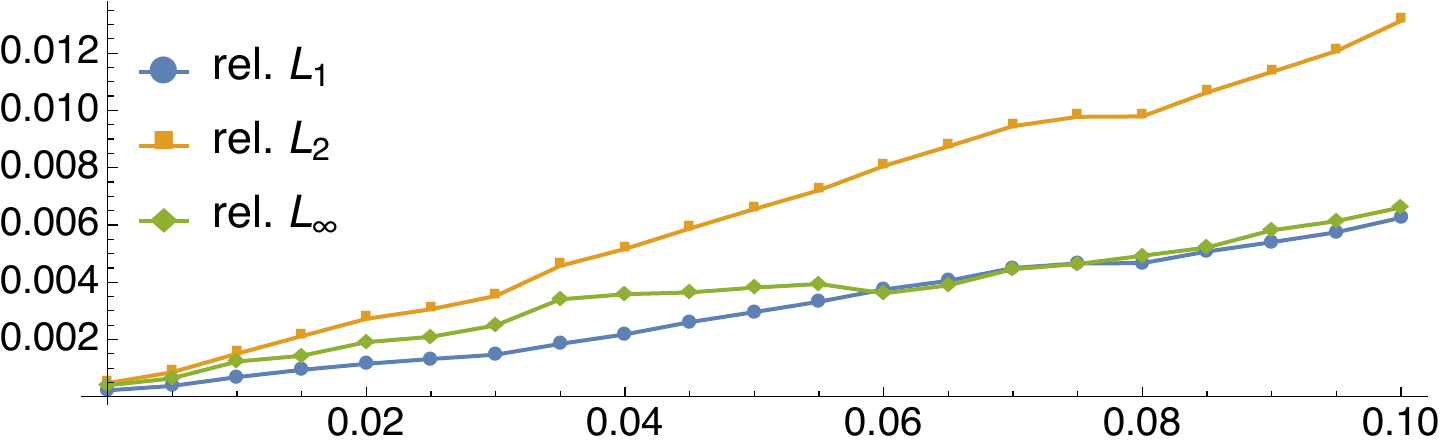}
  \includegraphics[width=0.475\textwidth]{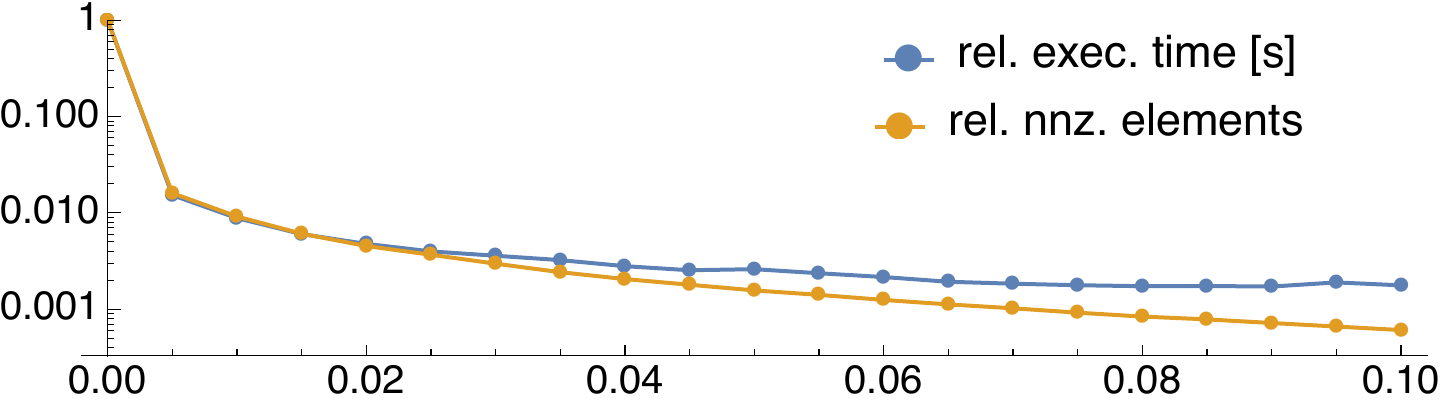}
  \caption{\emph{Top:} Reconstruction (blue) and $5X$-magnified error (rainbow colored) for the projection (black, dashed) of the signal in Fig.~\ref{fig:tomo:rec} for different hard thresholds. \emph{Bottom:} Corresponding relative error rates (left) with respect to the wavelet coefficients and execution time and non-zero elements (right) as a function of the threshold.}
  \label{fig:thresholded}
\end{figure}

\paragraph{Local projection}

As a basic proof-of-concept for our local Fourier slice equation we used a Gaussian for which an analytic solution for the projection is available, see Fig.~\ref{fig:gaussian}.
The errors compared to the reference are $L_1 = 3.78 \times 10^{-3}$, $L_2 = 2.0 \times 10^{-4}$, and $L_{\infty} = 1.51 \times 10^{-5}$.
The remaining residual results from the truncation of the basis to $[-10,10]^2$ and it can be decreased by increasing the apron region of the basis.

In Fig.~\ref{fig:gaussian} we also demonstrate the locality of our Fourier slice equation by reconstructing the Gaussian only over the positive $x_1$-axis.
Up to a small apron, the reconstruction only involves coefficients where $k_1 \geq 0$ so that the computational costs are $51.97\%$ of those for the full axis.
This is also reflected in the computation time that is $55\%$ of the full one, as expected from the theoretical analysis of Sec.~\ref{sec:construction:costs}.

\paragraph{Projection along arbitrary axes}

The foregoing example was isotropic, i.e. the projection along any direction yielded the same result.
As a simple non-isotropic example we consider the indicator function $\chi_S(x) \equiv \chi_{[-2.5,2.5]^2}(x)$ of a square.
Fig.~\ref{fig:square} shows the error as a function of the projection direction.
It can be seen that the error in the projected signal fluctuates but remains overall similar.
Interestingly, the smallest error is not attained for an axis-aligned projection but when $\theta_{\nu} = 45^{\circ}$.
The projected signal is then a tent function, which has a higher regularity than the box function obtained in the axis-aligned case, which explains the observed results.

\paragraph{Robustness of projection to $\epsilon$-thresholding}

Fig.~\ref{fig:thresholded} shows the projection error when small coefficients are (hard) thresholded to zero to increase the sparsity of the signal representation, as used for example in lossy image compression.
We see that thresholding can be exploited and that the projection error increases only moderately with the threshold.
For instance, with just $1.5\%$ nonzero coefficients, which, when stored as a sparse matrix, corresponds to about $4.1\%$ of the original memory requirements, we obtain an $L_{\infty}$ error that differs only by a factor of $1.6$ from the one without thresholding.
Fig.~\ref{fig:thresholded} also shows that the reduction in memory requirements translates directly into a reduction of the computation time.

The relative sparsity data reported in Fig.~\ref{fig:thresholded} is with respect to the full, quite redundant frame representation.
However, also compared to the uncompressed signal representation we only require $1/3$ of the original storage when $5\%$ of the coefficients are nonzero.

The Shepp-Logan-like signal with its cartoon-like structure is almost ideally suited for the curvelet-like polar wavelets used in the experiments, cf. Appendix~\ref{sec:appendix:implementation}.
For other signals, the sparse representation will require more coefficients to attain an acceptable error.
However, tomographic data, which is one of the principal applications of the Fourier slice theorem, also has a cartoon-like structure.

\subsection{Tomographic reconstruction}
\label{sec:experiments:tomography}

\begin{figure}
  \centering
  \includegraphics[width=1.0\textwidth]{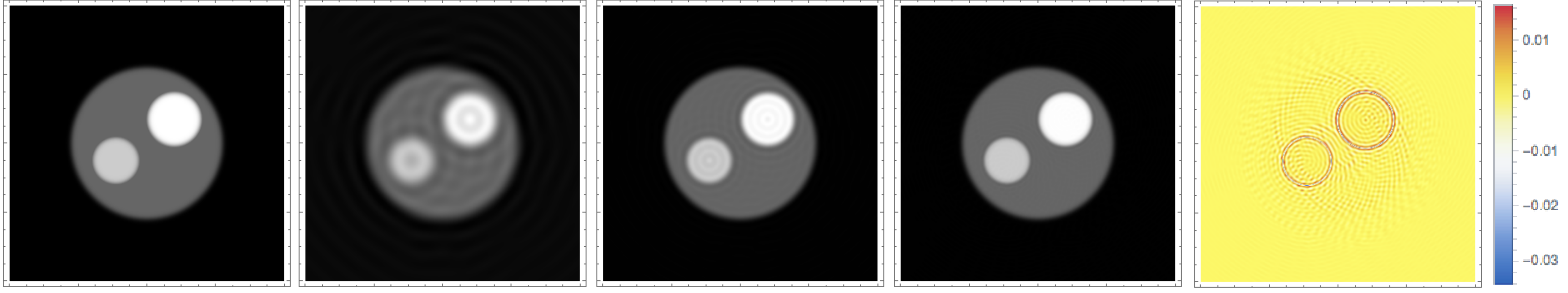}
  \caption{Tomographic reconstruction of the Shepp-Logan-like signal on the left for $j \leq 2,3,4$ levels (from left to right) with $196 \times 256$ measurements and basis functions in $[-5,5]^2$. The right plot shows the relative error for $j=4$.}
  \label{fig:tomo:rec}
\end{figure}

Tomographic reconstruction is the cornerstone of many medical imaging techniques~\cite{Kak2001,Herman2009,Natterer2001} and it plays an important role also in many other fields, e.g.~\cite{Salvo2003}.
In the following, we demonstrate how the local Fourier slice theorem can be used for (local) tomographic reconstruction and how it enables one to exploit sparsity.
Note that our results are only meant to be a proof-of-concept for the simplest setup and approach possible.
More work will be required to fully analyze the behavior, e.g. in the presence of noise, and to compete with state-of-the-art methods that have been optimized over the years.

Let $\varrho(x) : V \to \mathbb{R}$ be the density in an $n$-dimensional volume $V$ that is to be determined.
A tomographic measurement $m_{\nu}(y)$ is a $n-1$ dimensional signal given by
\begin{align}
  \label{eq:tomography:measurement}
  m_{\nu}(y)
  = \log{ \! \left( \frac{I_{\mathrm{in}}}{I_{\mathrm{out}}}\right) }
  = \int_{\mathbb{R}_{\nu}} \varrho(x) \, dx^n
\end{align}
where $\mathbb{R}_{\nu}$ is the normal line bundle over the domain of $m(y)$ which, for simplicity, we will assume to be Euclidean, and $I_{\mathrm{in}}$ and $I_{\mathrm{out}}$ are the emitted and received intensities, respectively.

Applying the local Fourier slice equation to Eq.~\ref{eq:tomography:measurement} we can write the measurements $m_{\nu}(y)$ as
\begin{align}
  \label{eq:tomography:slice}
  m_{\nu}(y) = \sum_{s \in \mathcal{I}} \varrho_s \, \psi_s^{n-1,\nu}(y)
\end{align}
where the $\varrho_s$ are the coefficients for the $n$-dimensional density $\varrho(x)$, which is hence given by
\begin{align}
  \label{eq:tomography:rho:rep}
  \varrho(x) = \sum_{s \in \mathcal{I}} \varrho_s \, \psi_s^n(x) .
\end{align}

\begin{figure}
  \centering
  \includegraphics[width=0.9\textwidth]{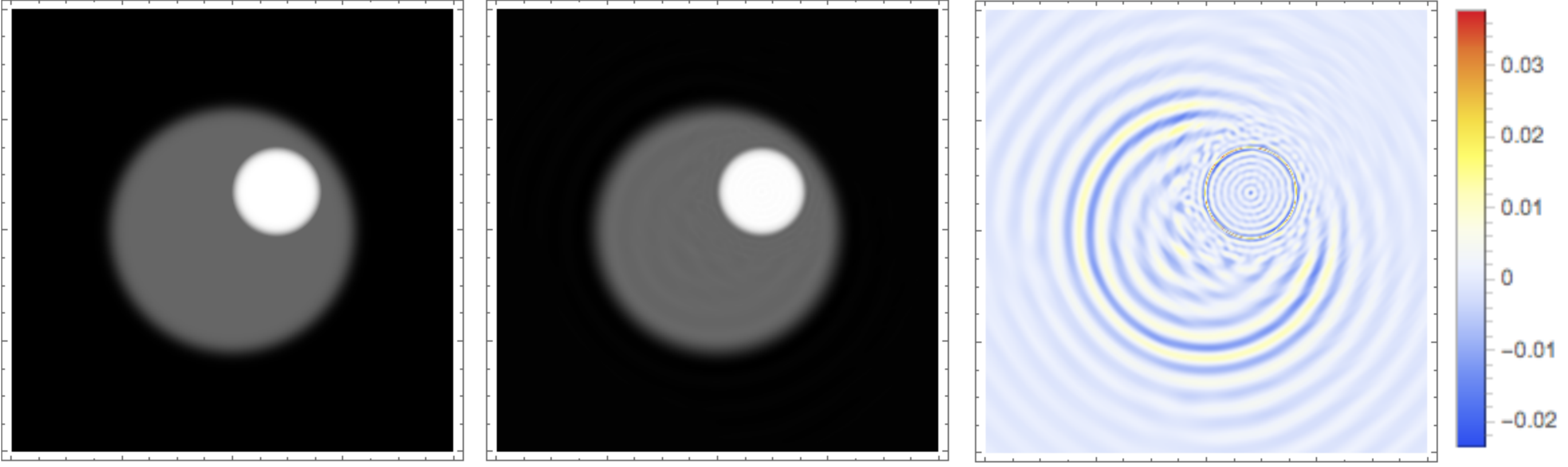}
  \caption{Sparse tomographic reconstruction (middle) of the Shepp-Logan-like signal (left) using the local Fourier slice equation for a sparse set of basis functions concentrated around the features, see Fig.~\ref{fig:tomo:sparse:sparsity}. The right plot shows the relative $L_{\infty}$ error. }
  \label{fig:tomo:sparse}
\end{figure}

Eq.~\ref{eq:tomography:slice} connects the frame coefficients $\varrho_s$, which specify the density we seek to reconstruct, with the projected signal $m_{\nu}(y)$, accessible through measurements.
We assume the measurements are pointwise values of $m_{\nu}(y)$ at locations $\lambda_i^{\nu} \in \Lambda_{\nu}$, i.e. $m_{\nu}(\lambda_i^{\nu})$.
Thus, with multiple measurement orientations $\nu_a$ we obtain a linear system whose rows are given by
\begin{align}
  \label{eq:tomography:ls}
  m_{\nu}(\lambda_{i}^{a_j}) = \sum_{s \in \mathcal{I}}  \varrho_s \, \psi_s^{n-1,\nu_a}\big( \lambda_{i}^{\nu_j} \big) .
\end{align}
In matrix-vector notation it takes the form $\mathrm{m} = \mathrm{Z} \, \mathrm{r}$ where $\mathrm{m}$ is the vector of all measured values $m_{\nu}(\lambda_{i}^{{\nu}_a})$, $\mathrm{Z}$ the matrix formed by the $\psi_s^{n-1,\nu_a}(\lambda_{i}^{\nu_j})$ with $i,j$ yielding the row index and $s$ the column one, and $\mathrm{r}$ the vector of the density coefficients $\varrho_s$ we seek to reconstruct.
In principle, the system $\mathrm{m} = \mathrm{Z} \, \mathrm{r}$ can be solved when the total number of measured values $\sum_a \vert \Lambda_{\nu_a} \vert$ satisfies $\sum_r \vert \Lambda_{a_r} \vert \geq \vert \mathcal{I} \vert$ and $\mathcal{I}$ is the index set in Eq.~\ref{eq:tomography:rho:rep}, although in practice successful reconstruction will require to use more measurements than unknowns.

\paragraph{Basic validation}

To validate the tomographic reconstruction in Eq.~\ref{eq:tomography:ls} we implemented it for the $2$-dimensional, Shepp-Logan-like test density $\varrho(x)$ shown in Fig.~\ref{fig:tomo:rec}.
Measurements were computed numerically using ray marching and for reconstruction we used a linear least squares fit (with a singular value cut-off of $10^{-6}$).
We used isotropic wavelets in $[-5,5]^2$ and up to level $j=4$ for this experiment so that the full basis consisted of $34,725$ functions.

Fig.~\ref{fig:tomo:rec} shows a level-by-level reconstruction.
When all levels $j=-1 ... 4$ are used then a visually artifact free image can be obtained and the remaining $L_{\infty}$ error is below $3\%$.
Some ringing occurs in the reconstruction when one considers the error plot.
But this is to be expect given that we use a linear least squared reconstruction that yields an $L_2$-projection suffering from the Gibbs phenomenon.

\begin{figure}
  \centering
  \includegraphics[width=0.86\textwidth]{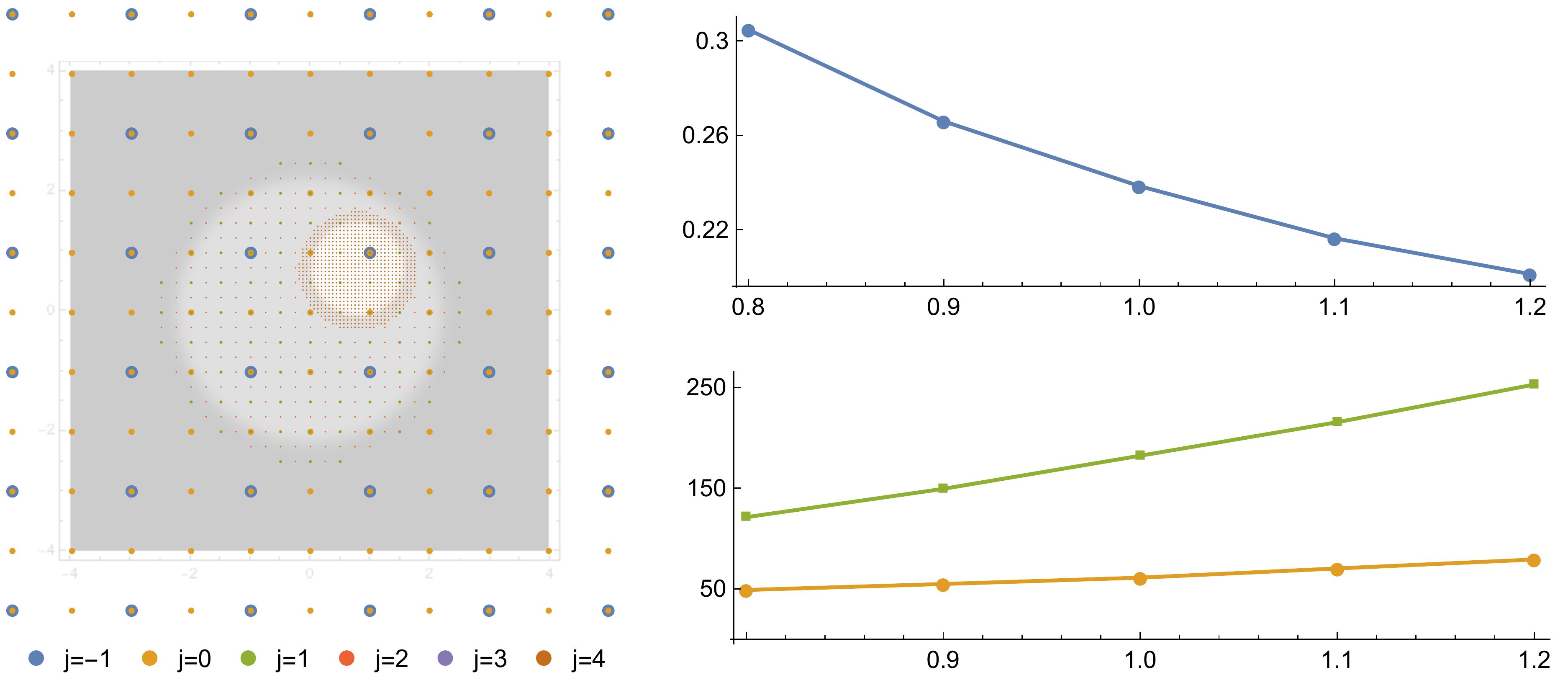}
  \caption{\emph{Left:} Sparse set of basis functions used for the reconstruction in Fig.~\ref{fig:tomo:sparse}. \emph{Right:} $L_2$ error (blue), computation time (yellow, in seconds), and memory (green, in MBs) as a function of the sparsity (relative, $1.0$ corresponds to the results in Fig.~\ref{fig:tomo:sparse}).}
  \label{fig:tomo:sparse:sparsity}
\end{figure}

\paragraph{Sparse reconstruction}
Fig.~\ref{fig:tomo:sparse} demonstrates sparse tomographic reconstruction where we exploit a priori regularity information  about the density signal.
The sparse set of basis functions we used in the experiment is shown in Fig.~\ref{fig:tomo:sparse:sparsity}, left.
On the scaling function level, $j=-1$, and the first wavelet level, $j=0$, we used all basis functions in $[-5,5]^2$.
On the following two levels, $j=1$ and $j=2$, only basis functions in the central region around the larger ball were used.
On the finest two levels, $j=3$ and $j=4$, we have basis functions only around the small ball.
The basis functions in the sparse representation were determined based on proximity to the two structures in the signal, i.e. we rely here on the known results on the sparsity of signal representations in wavelets~\cite{Meyer1992,Daubechies1992,Mallat2009}.
This resulted in $1588$ basis functions, compared to $34,725$ for the full basis.
For reconstruction we used $56$ orientations with $256$ samples.
Other parameters were as in the foregoing experiments.

The results in Fig.~\ref{fig:tomo:sparse} demonstrate that a faithful reconstruction of the density can be obtained using a priori sparsity information.
Compared to the full reconstruction in Fig.~\ref{fig:tomo:rec}, the sparse implementation requires only $1.4\%$ of memory and $2.1\%$ of computation time.
Some small ringing artifacts can be observed off the central region of interest where we focused our computations but the relative $L_{\infty}$ error remains below $3\%$.
It can be reduced by giving up some of the sparsity, as shown in Fig.~\ref{fig:tomo:sparse:sparsity}, right, where we plot the dependence of the reconstruction error, memory requirements and computation time on the sparsity, demonstrating that an effective trade-off between memory / time and reconstruction quality is possible.
A sparse basis suggests to also use a sparse set of samples that is also adapted to the input signal.
We leave this to future work.

The just presented approach for sparse reconstruction relies on a priori information about the signal.
This is not always available.
However, in applications such as medical imaging much is already known about the signal.
The alternative is to follow the methodology of compressed sensing and perform sparsity optimization, as in previous work on tomography using curvelets and shearlets~\cite{Frikel2013,Garduno2011,Vandeghinste2013}.
Such nonlinear optimization are, however, orders of magnitude more expensive than our approach.
The adaptive resolution used in the experiment can also be seen as a "magnifying lens" that enables one to reconstruct a high resolution signal in a region of interest.

\section{Conclusion}
\label{sec:conclusion}

We presented a local Fourier slice equation that, among other things, enables local projection and the effective use of sparsity.
Our construction exploits that wavelets defined in polar coordinates in frequency space respects the geometry of the Fourier slice theorem, since a slice is an iso-parameter set there.
With such wavelets one hence obtains a sequence closed under projection in that ``slicing'' a three-dimensional polar wavelet yields a two-dimensional one, and ``slicing'' a two-dimensional one yields a one-dimensional.
Moreover, all wavelets have closed form expressions in space and frequency, greatly facilitating their implementation.
This distinguishes polar wavelets and makes them a natural choice for our work.

In contrast to the classical Fourier slice theorem, our result does not require a discretization but is directly amenable to computations.
Furthermore, the computational complexity depends linearly on the size of the projection and the sparsity of the signal with respect to the projection direction, and scales hence in a natural manner with the complexity of the projected signal.
We demonstrated the validity of our local Fourier slice equation with experiments, in particular verifying its ability to project a signal only locally and to exploit sparsity to improve efficiency.
To show the potential of our theoretical result for practical applications, we considered tomographic reconstruction and  demonstrated that our local slice equation enables one to exploit sparsity, with substantial savings in memory and computation time.

Our work opens up many avenues for future work.
It would be interesting to obtain quantitative results on the conservation of sparsity that was experimentally demonstrated in Sec.~\ref{sec:experiments:basic}, e.g. based on the results on Quinto~\cite{Quinto1993,Quinto2007} that characterize the effect of the Radon transform on signal singularities.
As another application for our Fourier slice equation we envision image reconstruction of light field data following the approach proposed by Ng~\cite{Ng2005}, where the huge storage requirements make it essential that this can be performed directly from a compressed representation of the data.

Our experiments on tomographic reconstruction indicate that our local Fourier slice equation might be useful for this application.
However, much work remains to thoroughly evaluate its potential, in particular given that previous work that used harmonic bases, such as Haar wavelets and shearlets, ultimately did not improve over classical approaches~\cite{Garduno2011,Garduno2017}.
In our current implementation we did not use directional wavelets since these lead to a significantly larger redundancy.
Nonetheless, since they are very well suited for the signals that occur in medical applications of tomography, their use should be explored.
Iterative approaches are known to improve the accuracy of tomographic reconstruction~\cite{Beister2012}.
They are hence another obvious direction for future work.
An alternative would be to use operator-theoretic approaches, similar to~\cite{Candes2000,DeHoop2009a}.
J{\o}rgensen and co-workers~\cite{Jorgensen2015a,Jorgensen2015} studied sparse tomographic reconstruction and used concepts from compressed sensing~\cite{Candes2006a,Donoho2006a} to determine the number of required measurements.
It would be interesting to connect their results to our sparse reconstruction.
These authors also proposed a test set for sparse tomographic reconstruction which our approach should be evaluated on.
This set also contains noisy images, which is an important aspect we ignored.

%


\bibliography{fourierslice}

\appendix

\section{Conventions}
\label{sec:conventions}

The unitary Fourier transform of a function $f \in L_1(\mathbb{R}^n) \cap L_2(\mathbb{R}^n)$ is defined as
\begin{align}
\label{eq:fourier_transform}
  \mathcal{F}(f)(\xi) = \hat{f}(\xi) = \frac{1}{(2\pi)^{n/2}} \int_{\mathbb{R}_x^n} f(x) \, e^{-i \langle x , \xi \rangle} \, d x .
\end{align}
Its analogue on the sphere is the spherical harmonics expansion.
For $f \in L_2(S^2)$ it is given by
\begin{subequations}
  \label{eq:spherical_harmonics:expansion}
\begin{align}
  f(\omega)
  = \sum_{l=0}^{\infty} \sum_{m=-l}^l \langle f(\eta) , y_{lm}(\eta) \rangle \, y_{lm}(\omega)
  = \sum_{l=0}^{\infty} \sum_{m=-l}^l f_{lm} \, y_{lm}(\omega)
\end{align}
\end{subequations}
where $\langle \cdot , \cdot \rangle$ denotes the standard $L_2$ inner product on $S^2$.
We use standard (geographic) spherical coordinates with $\theta \in [0,\pi]$ being the polar angle and $\phi \in [0,2\pi]$ the azimuthal one.
The spherical harmonics basis functions $y_{lm}(\omega)$ in Eq.~\ref{eq:spherical_harmonics:expansion} are given by
\begin{align}
  \label{eq:spherical_harmonics}
  y_{lm}(\omega) = y_{lm}(\theta,\phi) = C_{lm} \, P_l^m( \cos{\theta} ) \, e^{i  m  \phi}
\end{align}
where the $P_l^m(\cdot)$ are associated Legendre polynomials and $C_{lm}$ is a normalization constant so that the $y_{lm}(\omega)$ are orthonormal over the sphere.

\section{Derivation of Fourier Slice Equation in Three Dimensions}
\label{sec:appendix:3d}

\paragraph{Projection along one dimension}

In three dimensions, the Fourier slice theorem for integration along the $x_3$-axis is given by
\begin{subequations}
  \label{eq:slice:3d:axis}
\begin{align}
  \label{eq:slice:3d:axis:1}
  f^{12}(x_{12})
  &= \frac{1}{(2\pi)^{3/2}} \int_{\mathbb{R}_{\xi_{12}}^2} \hat{f}(\xi_{12},0) \, e^{i \langle \xi_{12} , x_{12} \rangle} \, d\xi_{12}
\end{align}
where $\xi_{12} = (\xi_1,\xi_2,0)$.
Writing  $\hat{f}(\xi)$ using its polar wavelet representation in Eq.~\ref{eq:slice} we obtain
\begin{align}
  f^{12}(x_{12})
  &= \frac{1}{(2\pi)^{3/2}} \sum_{s \in \mathcal{I}} f_s \, \int_{\mathbb{R}_{\xi_{12}}^2} \hat{\psi}_s\big( x_{12},0 \big) \, e^{i \langle \xi_{12} , x_{12} \rangle} \, d\xi_{12}  .
\end{align}
The three-dimensional polar wavelet is given by Eq.~\ref{eq:polarlets:3d:hat} and using also the definition of the spherical harmonics in Eq.~\ref{eq:spherical_harmonics} we obtain
\begin{align}
  \hat{\psi}_s^{12}\big( \xi_{12},0 \big)
  &= \sum_{l,m} \, C_{lm} \, \kappa_{lm}^{jt} \, P_{lm}\big(\pi / 2\big) \, e^{i \, m \, \phi_{\xi_{12}}} \, \hat{h}\big(2^{-j} \vert \xi_{12} \vert \big) \, e^{-i \langle \xi_{12} , k_{12} \rangle} .
\end{align}
Re-arranging terms yiels
\begin{align}
  \hat{\psi}_s\big( \xi_{12},0 \big)
  &= \sum_{m} \underbrace{\left( \sum_l C_{lm} \, \kappa_{lm}^{jt} \, P_{lm}\big(\pi / 2\big) \right)}_{\displaystyle \beta_{m}^{jt}} \, e^{i \, m \, \phi_{\xi_{12}}} \, \hat{h}\big(2^{-j} \vert \xi_{12} \vert \big) \, e^{-i \langle \xi_{12} , k_{12} \rangle}
  \\[4pt]
  &= \sum_{m} \beta_{m}^{jt} \, e^{i \, m \, \phi_{\xi_{12}}} \, \hat{h}\big(2^{-j} \vert \xi_{12} \vert \big) \, e^{-i \langle \xi_{12} , k_{12} \rangle} .
\end{align}
\end{subequations}
For the general local Fourier slice equation, we can rotate the signal so that the projection direction is aligned with the $x_3$-axis, as in the calculation above, and then rotate the projected signal onto the plane $P_{\nu}$, which is trivial since this only pertains to where the projected signal is evaluated.
The rotated polar wavelet is given by
\begin{subequations}
\begin{align}
  \hat{\psi}_s^{\nu}( \xi_{\nu} )
  \equiv \hat{\psi}_s\big( R_{\nu} \xi_{12} \big)
  &= \hat{\gamma}\big( R_{\nu} \xi_{12} \big) \, \hat{h}(2^{-j_s} \vert \xi \vert) \, e^{i \langle R_{\nu} \xi_{12} , 2^{\textrm{--}j_s} k_s \rangle}
  \\[4pt]
  &=  \sum_{l,m} \left( \sum_{m'} W_{lm}^{m'}(\nu) \kappa_{lm'}^{j_s t_s} \right)  y_{lm}\big( \pi / 2 , \phi_{\xi} \big)
 \, \hat{h}(2^{-j_s} \vert \xi \vert) \, e^{i \langle R_{\nu} \xi_{12} , 2^{\textrm{--}j_s} k_s \rangle}
\end{align}
where in the second line we expressed the rotation in the spherical harmonics space where it is given by the elements of the Wigner-D matrices $W_{lm}^{m'}(\nu)$.
Writing the translation term as $ e^{i \langle R_{\nu} \xi_{12} , 2^{\textrm{--}j_s} k_s \rangle} = e^{i \langle \xi_{12} , 2^{\textrm{--}j_s} R_{\nu}^T k_s \rangle}$
\end{subequations}
we see that the rotation amount to working with a rotated grid $k_s^{\nu}= R_{\nu}^T k_s$ or, equivalently, with the projection of the unrotated grid $k_s$ onto $P_{\nu}$ spanned by $R_{\nu} \xi_{12}$.
Hence, the derivation of the axis-aligned case in Eq.~\ref{eq:slice:3d:axis} can be carried out unchanged using the rotated angular window coefficients $\kappa_{lm}^{j_s t_s,\nu}$ and for the rotated grid $k_s^{\nu}$.

\paragraph{Projection along two dimensions}

The Fourier slice theorem for the projection along two axes, say $x_1$ and $x_2$, is given by
\begin{subequations}
\begin{align}
  f^{3}(x_3)
  = \int_{R_{x_12}^2} f(x) \, d x_{12}
  = \frac{1}{(2\pi)^{3/2}} \int_{\mathbb{R}_{x_3}} \hat{f}(0,0,x_3) \, e^{i \langle x_3 , \xi_3 \rangle} \, d \xi_3 .
\end{align}
Inserting a polar wavelet representation for $\hat{f}(0,0,x_3)$ we obtain
\begin{align}
   f^{3}(x_3)
  = \frac{1}{(2\pi)^{3/2}} \sum_{s \in \mathcal{I}} f_s \int_{\mathbb{R}_{x_3}} \hat{\psi}_s(0,0,x_3) \, e^{i \langle x_3 , \xi_3 \rangle} \, d \xi_3 .
\end{align}
For the projected wavelet defined along the $x_3$-axis we have
\begin{align}
  \hat{\psi}_s^3(0,0,x_3)
  &= \underbrace{\left( \sum_{l} \kappa_{l}^{j_s t_s} \, y_{l}(0,0) \right)}_{\displaystyle \hat{\gamma}_{j_s,t_s}(0,0)} \, \hat{h}(2^{-j_s} \vert \xi \vert) \, e^{i \langle \xi_1 , 2^{\textrm{--}j_s} k_1^s \rangle}
\end{align}
where we used that the spherical harmonics satisfy $y_{lm}(0,0) = \delta_{l0}$.
Thus
\begin{align}
  \hat{\psi}_s^3(x_3)
  = \hat{\gamma}_{j_s,t_s}(0,0) \, \hat{h}(2^{-j_s} \vert \xi \vert) \, e^{i \langle \xi_1 , 2^{\textrm{--}j_s} k_1^s \rangle} .
\end{align}
\end{subequations}
The equivalent of the last equation for integration over an arbitrary plane follows from an argument analogous to those for the projection along one axis in three dimensions.

\subsection{Implementation details}
\label{sec:appendix:implementation}

In our implementation we use the radial window of the steerable pyramid~\cite{Portilla2000} since one then obtains closed form solutions for the radial profile $h_m(\vert x \vert)$ in the spatial domain, for the filter taps for the fast transform and other required quantities.
Since $h_m(\vert x \vert)$ is quite complicated we interpolate it at runtime.
For angular localization in $2D$ we used an extension of the wavelets for the interval proposed by Walter and Cai~\cite{Walter1999}, as first used in~\cite{Lessig2018a}, and with a curvelet-like coupling between the radial and angular scales.
This yields $1$, $4$, $8$, $12$ and $16$ different orientations on levels $j=0$ to $j=4$, respectively.
The scaling functions are always isotropic.
In $3D$ we used the spherical wavelets by McEwen~\cite{McEwen2016} and co-workers for angular localization.
We refer to the implementation in the supplementary material for further details.

\end{document}